\documentclass[10pt]{article}
\usepackage[francais,english]{babel}
\usepackage{amsmath}
\usepackage{amsfonts}
\usepackage{amssymb}
\usepackage{amsthm}
\usepackage{graphics}
\usepackage{amscd}
\usepackage{fullpage}
\usepackage[all]{xy}
\usepackage{epsfig}
\usepackage{color}
\newcommand{\R}{\mathbb{R}}

\newcommand{\U}{\mathbb{U}}
\newcommand{\E}{\mathbb{E}}

\newcommand{\mc}{\mathcal}

\newcommand{\mf}{\mathfrak}
\newcommand{\eps}{\varepsilon}
\newcommand{\ind}{{\bf 1}}
\renewcommand{\P}{\mathbb{P}}

\renewcommand{\H}{\mathbb{H}}

\DeclareMathOperator{\Id}{Id}

\DeclareMathOperator{\SLE}{SLE}

\title{Commutation relations for SLE}
\author{Julien Dub\'edat\footnote{Courant Institute}}
\newtheorem{thm}{Theorem}

\newtheorem{Thm}[thm]{Theorem}
\newtheorem{Def}[thm]{Definition}
\newtheorem{Prop}[thm]{Proposition}
\newtheorem{Lem}[thm]{Lemma}
\newtheorem{Cor}[thm]{Corollary}
\newtheorem{Conj}[thm]{Conjecture}

\begin{document}
\maketitle

\begin{abstract}
Schramm-Loewner Evolutions (SLEs) describe a one-parameter family of growth processes in the
plane that have particular conformal invariance properties. For
instance, SLE can define simple random curves in a simply connected
domain. In this paper we are interested in questions pertaining to the
definition of several SLEs in a domain (i.e. several random curves). 
In particular, one derives infinitesimal
commutation conditions, discuss some elementary solutions, study 
integrability conditions following from commutation
and show how to lift
these infinitesimal relations to global relations in simple cases. The
situation in multiply-connected domains is also discussed.
\end{abstract}

For plane critical models of statistical physics, such as percolation
or the Ising model, the general line of thinking of Conformal Field
Theory leads to expect the
existence of a non-degenerate scaling limit that satisfies conformal
invariance properties. Though, it is not quite clear how to define
this scaling limit and what conformal invariance exactly means.

One way to proceed is to consider a model in a, say, bounded (plane) simply-
connected domain with Jordan boundary, and to set boundary conditions
so as to force the existence of a macroscopic interface connecting two
marked points on the boundary. In this set-up, Schramm has shown that 
the possible scaling limits satifying conformal invariance along with
a ``domain Markov'' property are classified by a single positive
parameter $\kappa>0$, in the seminal article \cite{S0}. This defines the family of Schramm-Loewner
Evolutions (SLEs), that are probability measures supported on
non-self-traversing curves connecting two marked boundary points in a 
simply-connected domain.

Consider the following situation for critical site percolation on the
triangular lattice: a portion of the triangular lattice with mesh
$\eps$ approximates a fixed simply connected domain $D$ with two
points $x$ and $y$ marked on the boundary. The boundary arc $(xy)$ is
set to blue and $(yx)$ is set to yellow; sites are blue or yellow with
probability $1/2$. Then the interface between blue sites connected to
$(xy)$ and yellow sites connected to $(yx)$ is a non-self traversing
curve from $x$ to $y$. In this set-up, Smirnov has proved that the
interface converges to $\SLE_6$, as conjectured earlier by Schramm (\cite{Sm1,CamNew}).

For discrete models such as percolation or the Ising model, the full
information can be encoded as a collection of contours (interfaces
between blue any yellow, $+$ and $-$ spins, \dots). Hence it is
quite natural to consider scaling limits as collection of countours,
as in \cite{AizBur,CamNew}. Comparing the ideas of isolating one
macroscopic interface by setting appropriate boundary conditions
(following Schramm), and considering the scaling limit as a collection
of contours, one is led to the problem of describing the joint law of
a finite number of macroscopic interfaces created by appropriate
boundary conditions.

For instance, for percolation, consider a simply connected domain with
$2n$ marked points on the boundary, the $2n$ boundary arcs being
alternatively blue and yellow. This gives $n$ interfaces pairing the
$2n$ points. One can then consider the joint scaling limits of these
interfaces (either unconditionally or conditionally on a given
pairing). Each of these interfaces close to its starting point is
absolutely continuous w.r.t. $\SLE_6$. So we are defining $n$ ``non-crossing''
$\SLE_6$'s; the problem is then to precisely quantify their
interaction.

One remarkable feature of Schramm's construction is the classification
by a single positive parameter $\kappa$ for one interface satisfying
simple axioms. It is not hard to see that for, say, $2$ interfaces
connecting $4$ points, each interface is a priori described by
$\kappa$ and a drift term materializing the interaction. This drift
term can be seen as a function of the cross-ratio of the four boundary
points. 

The main goal of this article is to elucidate the constraints on the
drift terms imposed by the general geometric framework, and to prove
that in the most natural cases, the possible probability laws are
characterized by a finite number of parameters. The geometric
condition is that one can grow the interfaces in any order, at any
relative speed, and get the same result in distribution.

We aim at defining several SLEs in the same simply-connected
domain. As the growth of each SLE pertubates the time scales of other
SLEs, we want the collection of SLE to be invariant in distribution
under a global time reparametrization (that is,
a time change $\R_+^n\rightarrow \R_+^n$). At an infinitesimal level,
this invariance is expressed as a commutation relation for
differential operators (the infinitesimal generators of the driving
processes of the SLEs).

The conditions on the drift terms are non-linear differential
equations involving the drift terms pairwise. If one writes the drift
terms as log derivatives (in Girsanov fashion), then these conditions
can be written as a system of linear PDEs of a certain form satisfied by a single
``partition function''. The case where $2n$ points are marked on
the boundary and $n$ SLEs are grown is of particular interest.
The system of PDEs satisfied by the partition function is then unique and the
solution space is finite dimensional. The study of explicit solutions
is the subject of the companion paper \cite{D7}.

There is an analytically simple example of $n$ SLEs in a domain with
$(n+1)$ marked points. Here the drift terms are rational (and SLEs are
an example of $\SLE_\kappa(\underline \rho)$). In this case, we prove
that infinitesimal commutation conditions can be lifted to global
commutation; restriction formulae are derived for these. The radial
analogue is also discussed.

For multiply connected domains, as for simply-connected domains with,
say, more than 3 points marked on the boundary, the moduli space is no
longer a point, so identifying the drift terms that give ``physically
relevant'' SLEs is a problem. We work out the conditions imposed by
the following intuitive criterion: the interface can be grown ``from
both ends''. This involves ``cocycles'' on the configuration space.
The connection between this commutation condition and the
restriction property framework (\cite{LSW3,W3}) is made explicit.

These questions are connected with work by
Oded Schramm and David Wilson, Greg Lawler, and John Cardy; and also
by Bauer and Bernard for a more physical interpretation.

This paper is organized as follows. First, we review several examples
of commutation relations arising from different properties (general or
particular) of SLE. We then consider commutation at the level of
infinitesimal generators. Some (rational) solutions are discussed, and these
necessary infinitesimal conditions are recast as an integrability
problem. Study of integrability conditions leads to the classification
of some ``cocyles'' on the configuration space.
This leads to a holonomic system of PDEs (extending a natural
situation in critical percolation), when the residual moduli space is
reduced to a point. Restriction formulae are derived
in a particular (chordal) case; the analogue radial case follows. Finally, one discusses corresponding
questions in multiply connected domains.

\section{Introduction and notations}

First we recall some definitions and fix notations. We shall be mainly
interested in two kinds of SLE: chordal SLE in the upper half-plane
$\H$, from a real point to $\infty$; and radial SLE in the unit disk
$\U$, from a boundary point to $0$. Corresponding SLEs in other
(simply connected) domains are obtained by conformal equivalence. For
general background on SLE, see \cite{RS01,W1,Law}. Also, we will use
freely results on the restriction property and the ``loop soup''
(see \cite{LSW3,LW,W3}).

Consider the family of ODEs, indexed by $z$ in $\H$:
$$\partial_tg_t(z)=\frac 2{g_t(z)-W_t}$$ 
with initial conditions $g_0(z)=z$, where $W_t$ is some real-valued
(continuous) function. These chordal Loewner equations are defined up
to explosion time $\tau_z$ (maybe infinite). Define:
$$K_t=\overline{\{z\in\H:\tau_z<t\}}.$$
Then $(K_t)_{t\geq 0}$ is an increasing family of compact subsets of
$\overline\H$; moreover, $g_t$ is the unique conformal equivalence
$\H\setminus K_t\rightarrow \H$ such that (hydrodynamic normalization
at $\infty$):
$$g_t(z)=z+o(1).$$ 
For any compact subset $K$ of $\overline{\H}$ such that $\H\setminus
K$ is simply connected, we denote by $\phi$ the unique conformal
equivalence $\H\rightarrow \H\setminus K$ with hydrodynamic normalization
at $\infty$; so that $g_t=\phi_{K_t}$.

The coefficient of $1/z$ in the Laurent expansion of $g_t$ at $\infty$
is by definition the half-plane capacity of $K_t$ at infinity; this
capacity equals $(2t)$.

If $W_t=x+\sqrt\kappa B_t$ where $(B_t)$ is a standard Brownian
motion, then the Loewner chain $(K_t)$ (or the family $(g_t)$) defines
the chordal Schramm-Loewner Evolution with parameter $\kappa$ in
$(\H,x,\infty)$. The chain $K_t$ is generated by the trace $\gamma$, a
continuous process taking values in $\overline\H$, in the following
sense: $\H\setminus K_t$ is the unbounded connected component of
$\H\setminus\gamma_{[0,t]}$.

The trace is a continuous non self-traversing curve. It is a.s. simple
if $\kappa\leq 4$ and a.s. space-filling if $\kappa\geq 8$.

In the radial case, Loewner's equations are indexed by $z\in\U$,
$$\partial_tg_t(z)=-g_t(z)\frac {g_t(z)+\xi_t}{g_t(z)-\xi_t}$$
and $g_0(z)=z$, $\xi$ takes values in the unit circle. The hull $K_t$
is defined as above, and $g_t$ is the unique conformal equivalence
$\U\setminus K_t\rightarrow \U$ with $g_t(0)=0$,
$g'_t(0)>0$. Moreover, $g'_t(0)=e^{-t}$. If
$\xi_t=\xi_0\exp(i\sqrt\kappa B_t)$, where $B$ is a standard Brownian
motion, one gets radial $\SLE_\kappa$ from $\xi_0$ to $0$ in $\U$.

Note that chordal SLE depends only on two boundary points, and radial
SLE depends on one boundary and one bulk point. In several natural
instances, one needs to track additional points on the boundary. This
has prompted the introduction of $\SLE(\kappa,\rho)$ processes in
\cite{LSW3}, generalized in \cite{D4}. The driving Brownian motion is
replaced by a semimartingale which has local Girsanov density w.r.t.
the original Brownian motion.

In the chordal case, let $\underline\rho$ be a multi-index, i.e. :
$$\underline\rho\in \bigcup_{i \ge 0} \R^i$$
Let $k$ be the length of $\underline\rho$; if $k=0$, one simply
defines $\SLE(\kappa,\varnothing)$ as a standard $\SLE_\kappa$. If
$k>0$, assume the existence of processes $(W_t)_{t\geq 0}$ and
$(Z^{(i)}_t)_{t\geq 0}$, $i\in\{1\dots k\}$ satisfying the SDEs:
\begin{equation}\label{E1}
\left\{\begin{array}{l}dW_t=\sqrt\kappa dB_t+\sum_{i=1}^k\frac{\rho_i}{W_t-Z^{(i)}_t}dt\\
dZ^{(i)}_t=\frac 2{Z^{(i)}_t-W_t}dt\end{array}\right.
\end{equation}
and such that the processes $(W_t-Z^{(i)}_t)$ do not change sign.
Then we define the chordal $\SLE_\kappa(\underline\rho)$ process starting from
$(w,z_1,\dots z_k)$ as a chordal Schramm-Loewner evolution the driving process
of which has the same law as $(W_t)$ as defined above, with
$W_0=w,Z^{(i)}_0=z_i$.

In the radial case, assume the existence of processes $(\xi_t)_{t\geq 0}$ and
$(\chi^{(i)}_t)_{t\geq 0}$, $i\in\{1\dots k\}$ satisfying the SDEs:

\begin{equation}\label{E2}
\left\{\begin{array}{l}d\xi_t=(i\xi_t\sqrt\kappa dB_t-\frac\kappa
2\xi_tdt)+\sum_{i=1}^k\frac{\rho_i}2\left(-\xi_t\frac{\xi_t+\chi^{(i)}_t}{\xi_t-\chi^{(i)}_t}\right)dt\\
d\chi^{(i)}_t=-\chi^{(i)}_t\frac{\chi^{(i)}_t+\xi_t}{\chi^{(i)}_t-\xi_t}dt\end{array}\right.
\end{equation}
The processes $\xi$, $\chi^{(i)}$ may bounce on each other but not cross. 
This defines radial $\SLE_\kappa(\underline\rho)$ in the unit
disk. Note the factor $1/2$ before the $\rho_i$ parameters in the SDE:
this is to ensure coherence with the chordal case.

\section{Examples of commutation}

We begin by discussing how properties of SLE (e.g. reversibility
and duality) yield natural examples of commutation relations.

{\bf Reversibility:} consider a chordal SLE in $(\H,0,\infty)$,
$\gamma$ its trace. For simplicity, assume that $\kappa\leq 4$, so
that the trace is a.s. simple. Define
$\hat\gamma_t=\gamma_{1/t}$. Then $\hat\gamma$ is a simple curve from
$\infty$ to $0$ in $\H$ (for transience of SLE, see
\cite{RS01}). After a time change $s=s(t)$, $\hat\gamma$ is such that
$-1/\hat\gamma_{[0,s]}$ has capacity $2s$. Then, according to
reversibility, $\hat\gamma$ a $\SLE_\kappa$ in $\H$ from $\infty$ to
$0$. (Note that for all $0<\kappa\leq 4$ and involution of $\H$ of
type $z\mapsto -\lambda/z$, this defines a somewhat intricate
measure-preserving involution of the Wiener space). 

Admitting reversibility, one can define a chordal SLE growing ``from
both ends'' in the following fashion: let $B$ be a standard Brownian
motion, with filtration ${\mc F}$, $\gamma$ the trace of the
associated $\SLE_\kappa$, and $\hat\gamma$ as above. Then consider
$(K_{t,s}=\gamma_{[0,t]}\cup\hat\gamma^{-1}_{[0,s]})_{s,t\geq 0}$ and
the filtration $({\mc G}_{t,s})$ it generates. For any $(t_0,s_0)$,
$(K_{t_0+t,s_0})_t$ is a (time-changed) chordal SLE in $\H\setminus
K_{t_0,s_0}$, from $\gamma_{t_0}$ to $\hat\gamma_{s_0}$; conversely 
$(K_{t_0,s_0+s})_s$ is a (time-changed) chordal SLE in $\H\setminus
K_{t_0,s_0}$, from $\hat\gamma_{s_0}$ to $\gamma_{t_0}$. 

Together with conformal equivalence, this gives the following Markov property: if $f_{t,s}$ is a conformal
equivalence $(\H\setminus K_{t,s},\gamma_t,\hat\gamma_ s)\rightarrow
(\H,0,\infty)$ with some normalization (e.g. $f_{t,s}(1)=1$), then
$f_{t_0,s_0}(K_{t_0+t,s_0+s})$ is up to a time-change
$\R_+^2\rightarrow \R_+^2$ a copy of $(K_{t,s})$ independent from
${\mc G}_{t_0,s_0}$.

{\bf Duality:} Duality relates the boundary of non-simple SLE
($\kappa>4$) with corresponding simple SLEs ($\kappa'=16/\kappa$). Let
us try to formulate a precise conjecture in a ``dual'' fashion. We
elaborate on restriction formulae identities discussed in \cite{D4}. 

Let $\kappa>4$, $\kappa'=16/\kappa$. Consider the configuration $(\H,x,y,z,\infty)$, where $x<y<z$. Define a Loewner chain from $0$ to $\infty$ as follows:

the chain $(K_t)_{t\leq\tau_z}$ is an
$\SLE_\kappa(\kappa/2-4,-\kappa/2)$ in $\H$, started from $(x,y,z)$,
aiming at $\infty$, and stopped at time $\tau_z$ when the trace hits
$z$ (which it does with probability 1). Then $(K_{t+\tau_z})_{t\geq
0}$ is a $\SLE_\kappa(\kappa-4)$ in $\H\setminus K_{\tau_z}$, started
from $(z,z^+)$ and aiming at $\infty$.

The right-boundary of $K_\infty=\bigcup_{t\geq 0}K_t$ is a simple
curve from $z$ to $\infty$ in $\H$; denote by $(\delta_u)$ the
corresponding Loewner trace (i.e. $\delta_{[0,\infty)}$ is the
right-boundary of $K_\infty$ and $\delta_{[0,u]}$ has half-plane
capacity $(2u)$\ ).

Now consider a configuration $(\H,x',y',z',\infty)$, where $x'<y'<z'$,
and let $\gamma'$ be the trace of the chordal
$\SLE_{\kappa'}(-\kappa'/2,\kappa'-2)$ in $\H$ started from
$(z',x',y')$ and aiming at $\infty$.

Then we can formulate:

\begin{Conj} The following statements hold:\\
(i) The law of $\delta$ is that of $\gamma$, where $x=x'$, $y=y'$, $z=z'$.\\
(ii) The law of $(\phi_{\delta_{[0,u]}}(K_t))_{t\geq 0}$ conditionally
on $\delta_{[0,u]}$ is (up to a time-change) that of a copy of $(K_t)$
started from $(x',y',z')=\phi_{\delta_{[0,u]}}(x,y,\delta_u)$.
\end{Conj}

This conjecture can be interpreted in terms of multiple SLEs: one can
grow simultaneously the chain $(K_t)$ and its (final)
right-boundary. One also get a Markov property similar to the one
discussed for reversibility.

{\bf Locality:} The scaling limit of the exploration process for
critical site percolation on the triangular lattice is $\SLE_6$
(see \cite{Sm1,CamNew}). For some boundary conditions, one can define several
exploration paths. Consider for instance the following situation:
$(D,x_1,\dots, x_{2n})$ is a simply connected domain with $(2n)$ marked
boundary points in cyclical order. The segments
$(x_1,x_2),\dots,(x_{2n-1},x_{2n})$
(resp. $(x_2,x_3),\dots,(x_{2n},x_{1})$) are set to blue
(resp. yellow). Then one can start an exploration process at each of
the boundary points $x_i$; these are well-defined up to some
disconnection event.

One can also consider some conditional versions: for instance,
critical percolation in $(\H,0,1)$, where the half-lines $(\infty,0)$ and
$(1,\infty)$ are blue and $(0,1)$ is yellow, conditionally on the
existence of a yellow path from $(0,1)$ to infinity (this is a
singular conditioning, related to the one-arm half-plane
exponent). Now the exploration processes started from $0$ and $1$
resp. can be defined for all time. The two traces intersect at pivotal
points for the conditioning event.
 
One may also consider the following
situation: a conformal rectangle, with sides alternately blue and
yellow. Hence, one can start four exploration processes (one at each
vertex). Then condition on a Cardy crossing event (e.g. the two blue
sides are connected by a blue path). One can note that in this
situation, the Girsanov drift terms are not rational functions.

{\bf Restriction:} The restriction property of $\SLE_{8/3}$ can be
used to get commutation relations. For instance, consider a simply
connected domain with four marked points on the boundary, say
$(\H,a,b,c,d)$. One can define two independent $\SLE_{8/3}$'s, from
$a$ to $b$ and $c$ to $d$ resp., and condition them on not
intersecting. Then, from the restriction property, this system of two
SLEs has a natural Markov property, and also a restriction property. 

More precisely, let $\gamma$ and $\gamma'$ be the traces of these
SLEs, $(g_t)$ the family of conformal equivalences of the first one
(for some time parameterization). Then (${\mc L}$ denotes
probability distributions)
\begin{align*}
{\mc L}(g_t(\gamma_{(t,\infty)}),g_t(\gamma')|\gamma\cap\gamma'=\varnothing)
&={\mc
L}(g_t(\gamma_{(t,\infty)}),g_t(\gamma')|\gamma_{(0,t)}\cap\gamma'=\varnothing,g_t(\gamma_{(t,\infty)})\cap
g_t(\gamma')
)\\
&={\mc
L}(\tilde\gamma,\tilde\gamma'|\tilde\gamma\cap
\tilde\gamma'=\varnothing)
\end{align*}
where $\tilde\gamma$ and $\tilde\gamma'$ are independent $\SLE_{8/3}$'s
going from $g_t(\gamma_t)$ to $g_t(b)$ and from $g_t(c)$ to $g_t(d)$
resp. (using the Markov property for $\gamma$ and the restriction
property for $\gamma'$). For the restriction property, note that, for
any hull $A$ disjoint from $\{a,b,c,d\}$:
\begin{align*}
\P((\gamma\cup\gamma)'\cap
A=\varnothing|\gamma\cap\gamma'=\varnothing)&=\frac{\P(\gamma\cap\gamma'=\gamma\cap
A=\gamma'\cap A=\varnothing)}{\P(\gamma\cap\gamma'=\varnothing)}\\
&=\frac{\P(\gamma\cap\gamma'=\varnothing|\gamma\cap
A=\gamma'\cap A=\varnothing)\P(\gamma\cap
A=\gamma'\cap A=\varnothing)}{\P(\gamma\cap\gamma'=\varnothing)}\\
&=\frac{\P(\phi_A(\gamma)\cap\phi_A(\gamma')=\varnothing)\P(\gamma\cap
A=\gamma'\cap A=\varnothing)}{\P(\gamma\cap\gamma'=\varnothing)}\\
&=\P(\gamma\cap
A=\varnothing)\P(\gamma'\cap A=\varnothing)\frac{\psi(\phi_A(a,b,c,d))}{\psi(a,b,c,d)}
\end{align*}
where $\psi(a,b,c,d)$ is the probability that the two independent
$\SLE_{8/3}$ do not intersect.

If $\kappa\in (0,8/3)$, one can consider two independent
$\SLE_\kappa$, a corresponding independent loop soup, and condition on the event:
no loop intersects the two SLEs. A standard computation shows that the
probability $\psi_\kappa$ of this event is given by:
$$\psi_\kappa(\infty,0,x,1)=\frac{\Gamma(4/\kappa)\Gamma(12/\kappa-1)}{\Gamma(8/\kappa)\Gamma(8/\kappa-1)}x^{2/\kappa}\hphantom{F}_2F_1\left(\frac 4\kappa,1-\frac
4\kappa;\frac 8\kappa;x\right)$$

In a domain with $(2n)$ marked points on the boundary in cyclical
order, say $(\H,x_1,\dots x_{2n})$, for a given pairing of $\{x_1,\dots
x_{2n}\}$, define $n$ independent $\SLE_\kappa$, with endpoints
determined by the pairing. Consider $(n-1)$ auxiliary independent loop
$L_1,\dots, L_{n-1}$
soups, with intensity $\lambda_\kappa$.
One can consider the event: for $k=1,\dots,n-1$, no loop in $L_k$
intersects more than $k$ SLEs. 
This has positive probability iff the pairing is a non-crossing
one. As is well known, there are $C_n$ of these pairings, where $C_n$
is the $n$-th Catalan's number:
$$C_n=\frac{{{2n}\choose{n}}}{n+1}$$
Then one can condition on this event to get $n$ non intersecting
$\SLE_\kappa$'s, that have an appropriate Markov property and
restriction property. This situation is discussed in details in
\cite{D7}; when $\kappa=2$, this is directly connected to
Fomin's formulae \cite{Fomin,KozL}.

{\bf Wilson's algorithm:} In the case $\kappa=8$, $\kappa'=2$, the
Uniform Spanning Tree and the Loop-Erased Random Walk converge to
$\SLE_8$ and $\SLE_2$ resp. (see \cite{LSW2}). As pointed out in
\cite{LSW2}, duality follows from these convergence and Wilson's
algorithm, that gives an exact relation between UST and LERW at the
discrete level (\cite{Wi}). 

Let us formulate a precise duality identity in this
situation. Consider $(K_t)$ a chordal $\SLE_8$ in $(\H,0,\infty)$. Let
$G$ be the (random) leftmost point visited by this $\SLE$ before
$\tau_1$. Then a standard SLE computation (see e.g. \cite{W1}) yields:
$$\P(G\leq -g)=\frac 1\pi\int_0^{g}\frac{dt}{(1+t)\sqrt t}=\frac
2\pi\arctan(\sqrt g)$$
This distribution is the exit distribution of a random walk with
normal reflection on $\R^+$, absorbed on $\R^-$, and started from 1 (as is
readily seen by mapping $\H$ to a quadrant by $z\mapsto\sqrt z$ and a
reflection argument).

At the discrete level, we are considering a UST wired on $\R^-$ and
free and $\R^+$. The branch connecting $1$ to $\R^-$ is a LERW started
from $1$ and reflected on $\R^+$. By a slight modification of the
arguments of \cite{LSW2} (considering the Poisson kernel for this
reflected random walk gives ``harmonic martingales'' for the
time-reverted LERW), one gets that conditionally on $G$, the boundary
of $K_{\tau_1}$, which is a random simple curve connecting $G$ and $1$,
and the scaling limit of this LERW, is chordal $\SLE_2(-1,-1)$ in $(\H,G,1)$
started from $(G,0,\infty)$.

Wilson's algorithm gives more information. The boundary $\partial
K_{\tau_1}$ divides $\H$ in two simply connected domains $\H_0$ and
$\H_\infty$, with $0$ and $\infty$ in their respective boundary. Then,
conditionally on $\partial
K_{\tau_1}$, the original $\SLE_8$ is the concatenation of a chordal $\SLE_8$
in $(\H_0,0,1)$ and a chordal $\SLE_8$
in $(\H_\infty,1,\infty)$.

For small times, the law of the original $\SLE_8$ conditionally on
$G=g$ (in the regular conditional probability sense) is easy to work
out. Consider the following martingale (with usual notations):
$$\arctan\sqrt\frac{W_t-g_t(x)}{g_t(1)-W_t}$$
Differentiating w.r.t $x$, one gets a local martingale:
$$\frac{g'_t(x)}{g_t(1)-g_t(x)}\sqrt\frac{W_t-g_t(x)}{g_t(1)-W_t}$$
Using this as a Girsanov density, one finds that the conditional
$\SLE_8$ is a chordal $\SLE_8(-4,4)$ in $(\H,0,\infty)$ started from $(0,G,1)$.

For symmetry, and from reversibility, the chordal $\SLE_8$
in $(\H_\infty,1,\infty)$ is a time-reversed  chordal $\SLE_8$
in $(\H_\infty,\infty,1)$.

Consider now these different processes as chordal SLEs in $\H$ aiming
at 1 (and not $\infty$). Then we have \\
an $\SLE_2(-1,-1)$ started from $(G,0,\infty)$\\
an $\SLE_8(-4,2)$ started from $(0,G,\infty)$\\
an $\SLE_8(-4,2)$ started from $(\infty,G,0)$\\
We shall see later that this fits in infinitesimal relations for $\SLE_\kappa(\underline\rho)$.

One can extend the situation as follows: in the discrete setting,
consider $n$ points on $\R^+$ and $n$ points on $\R^-$:
$$x_n<\cdots <x_1<0<y_1<\cdots<y_n.$$
Consider a UST with the same boundary conditions as before, and the
smallest subtree containing $y_1,\dots,y_n$ and $\R^-$. Condition on
the event that this subtree has no triple point in the bulk. Then it
consists in the union of $n$ disjoint paths in the bulk and $\R^-$. Now
condition on the endpoints of these branches being $x_n,\dots,x_1$,
and take this to the scaling limit. Using Wilson's algorithm and
Fomin's formulae (\cite{Fomin}), everything can be made explicit, and
this defines $n$ ``non-intersecting'' $\SLE_2$'s in the upper half-plane
(with $(2n+2)$ marked points on the boundary).

\section{Commutation of infinitesimal generators}

\subsection{The commutation framework}

We have seen natural examples where two SLEs could be grown in a
common domain in a consistent fashion. In this section, we discuss
necessary infinitesimal conditions. We shall define a ``global''
commutation condition, of geometric nature, and express its consequence
in terms of infinitesimal generators, which is of algebraic nature.

Let us consider the following chordal situation: the domain is $\H$, SLEs aim at
$\infty$, and $(x,y,z_1,\dots,z_n)$ are $(n+2)$ (distinct) points on
the real line; the point at infinity is also a marked point. We want to grow two infinitesimal hulls (with capacity of
order $\eps$) at $x$ and $z$ respectively. 
We can either grow a hull $K_\eps$ at $x$, and then another one at $y$
in the pertubed domain $\H\setminus K_\eps$, or proceed in the other
order. The coherence condition is that these two procedures yield the
same result.

Let us make things more rigorous. Consider a Loewner
chain $(K_{s,t})_{(s,t)\in{\mc T}}$ with a double time index, so that
$K_{s,t}\subset K_{s',t'}$ if $s'\geq s$, $t'\geq t$ and $K_{s,t}\neq
K_{s',t'}$ if $(s',t')\neq (s,t)$. We only consider chains up to time
reparameterization $\R_+^2\rightarrow\R_+^2$. We also assume that
$K_{s,t}=K_{s,0}\cup K_{0,t}$.
The time
set ${\mc T}$ may be random, but includes a.s. a neighbourhood of
$(0,0)$ in $\R_+^2$. Also, if $s\leq s'$, $t\leq t'$, $(s',t')\in{\mc
T}$, then $(s,t)\in{\mc T}$. Define
$g_{s,t}$ the conformal equivalence $\H\setminus K_{s,t}\rightarrow
\H$ with hydrodynamic normalization at infinity
($g_{s,t}=\phi_{K_{s,t}}$ with the earlier notation), and the continuous
traces $\gamma$, $\tilde\gamma$, such that:
$$\gamma_{s,t}=\lim_{\eps\searrow 0}\overline{ g_{s,t}(K_{s+\eps,t}\setminus
K_{s,t})}, \tilde\gamma_{s,t}=\lim_{\eps\searrow 0} \overline{g_{s,t}(K_{s,t+\eps}\setminus
K_{s,t})}$$
where  $\gamma_{0,t}=x$ for
all $(0,t)\in{\mc T}$, and similarly $\tilde\gamma_{s,0}=y$  for
all $(s,0)\in{\mc T}$. 
%Define also $K^x_{s,t}$ the compact subset of
%$\overline\H$ such that $\H\setminus K^x_{s,t}$ is the unbounded
%connected component of $\H\setminus \gamma_{[0,s]\times\{t\}}$ and
%$K^x_{s,t}\cap\R=\partial(K^x_{s,t}\cup\H)$; also
%$g^x_{s,t}=\phi_{K_{s,t}}$. The hulls $K^y_{s,t}$ and conformal
%equivalences $g^y_{s,t}$ are defined likewise.  

Furthermore, assume that 
the following conditions are satisfied:

\begin{Def}
Let $(K_{s,t})_{(s,t)\in{\mc T}}$ be a random Loewner chain with
double time indexing; the associated conformal equivalences are
$g_{s,t}=\phi_{K_{s,t}}$. We say that $(K_{s,t})$ is an
$\SLE(\kappa,b,\tilde\kappa,\tilde b)$ if:
\begin{enumerate}
\item The time set ${\mc T}$ is a.s. open, connected, and a
neighbourhood of ${(0,0)}$ in $\R_+^2$. The ranges of the traces $\gamma_{\mc T}$, $\tilde\gamma_{\mc
T}$ are disjoint and $z_1,\dots, z_n\notin K_{s,t}$ for $(s,t)\in{\mc T}$. 
\item Let $\sigma$ (resp. $\tau$) be a stopping time in the filtration
generated by $(K_{s,0})_{(s,0)\in{\mc T}}$
(resp. $(K_{0,t})_{(0,t)\in{\mc T}}$). Let also ${\mc T'}=\{(s,t):
(s+\sigma,t+\tau)\in {\mc T}\}$
and $(K'_{s,t})_{(s,t)\in{\mc
T'}}=\left(\overline{g_{\sigma,\tau}(K_{s+\sigma,t+\tau}\setminus
K_{s,t})}\right)$. Then $(K'_{s,0})_{(s,0)\in{\mc T'}}$ is distributed
as a stopped $\SLE_\kappa(b)$, i.e an $\SLE$ driven by:
$$dX_s=\sqrt{\kappa} dB_s+ b(X_s, g_s(y),\dots, g_s(z_i),\dots)dt$$ 
Likewise $(K'_{0,t})_{(0,t)\in{\mc T'}}$ is distributed
as a stopped $\SLE_{\tilde\kappa}(\tilde b)$, i.e an $\SLE$ driven by: 
$$dY_t=\sqrt{\tilde\kappa} d\tilde B_t+ \tilde b(\tilde g_t(
x),Y_t,\dots, \tilde g_t(z_i),\dots)dt$$

%\item The chain $(K_{s,0})_{(s,0)\in{\mc T}}$ is distributed as an $\SLE_\kappa(b)$, i.e. an $\SLE$ driven by
%$$dX_s=\sqrt{\kappa} dB_s+ b(X_s, g_s(y),\dots, g_s(z_i),\dots)dt$$ 
%stopped when the trace hits one of the marked points. Similarly, the
%chain  $(K_{0,t})_{(0,t)\in{\mc T}}$ is an $\SLE_{\tilde\kappa}(\tilde
%b)$, i.e. an $\SLE$ driven by
%$$dY_t=\sqrt{\tilde\kappa} d\tilde B_t+ \tilde b(\tilde g_t(
%x),Y_t,\dots, \tilde g_t(z_i),\dots)dt$$ 
%\item Let $\sigma$ be a stopping time in the filtration generated by
%$(K_{s,0})_{(s,0)\in{\mc T}}$. Let ${\mc
%T}_\sigma=\{(s,t):(s+\sigma,t)\in{\mc T}\}$. Then the chain
%$$\left(\overline {g_{\sigma,0}(K_{s+\sigma,t}\setminus
%K_{s,t})}\right)_{(s,t)\in{\mc T}_\sigma}$$
%is distributed as $(K_{s,t})_{(s,t)\in{\mc T}}$ (up to time
%reparameterization) and is independent of $K_{\sigma,0}$ conditionally
%on its initial state (postion of the marked points). Similarly, if $\tau$ is a stopping time in the filtration generated by
%$(K_{0,t})_{(0,t)\in{\mc T}}$, ${\mc
%T}_\tau=\{(s,t):(s,t+\tau)\in{\mc T}\}$, then the chain
%$$\left(\overline {g_{0,\tau}(K_{s,t+\tau}\setminus
%K_{s,t})}\right)_{(s,t)\in{\mc T}_\tau}$$
%is distributed as $(K_{s,t})_{(s,t)\in{\mc T}}$ (up to time
%reparameterization) and is independent of $K_{0,\tau}$ conditionally
%on its initial state.
\end{enumerate}
\end{Def}

Here $B$, $\tilde B$ are standard Brownian motions, 
$(g_s)$, $(\tilde g_t)$ are the associated conformal equivalences, $b$, $\tilde
b$ are some
smooth, translation invariant, and homogeneous of degree $(-1)$
functions. If 
$A^x$, $A^y$ are two increasing functions of hulls growing at $x$ and $y$
resp. (e.g. the half-plane capacity), we shall be particularly
interested in stopping times of type $\sigma=\inf(s:A^x(K_{s,0})\geq
a^x)$, $\tau=\inf(t:A^y(K_{0,t})\geq
a^y)$.

Note that $(X_s,\dots, g_t(z_i),\dots)$ is a Markov process. Let $P$ be its
semigroup and ${\mc L}$ its infinitesimal generator.
Similarly, $(\tilde g_t(\hat x),Y_t,\dots)$ is a Markov process
with semigroup $Q$ and infinitesimal generator ${\mc M}$. We are
interested in what conditions on the functions $b$ and $\tilde b$ are
implied by these assumptions (the existence of an $\SLE(\kappa,b,\tilde
\kappa,\tilde b)$).

So let $F$ be a test function $\R^{n+2}\rightarrow \R$, and $c>0$ be
some constant (ratio of speeds). We apply the previous assumptions
with $A^x=A^y={\rm cap}$ (the half-plane capacity), $a^x=2\eps$,
$a^y=2c\eps$. We are interested in the hull $K_{\sigma,\tau}$. Two ways
of getting from $K_{0,0}$ to $K_{\sigma,\tau}$ are (symbolically):
$$
\xymatrix{K_{0,0}\ar[r]\ar[d]& K_{\sigma,0}\ar[d]\\
K_{0,\tau}\ar[r]&K_{\sigma,\tau}
}$$
and our assumptions give a description of these transitions.

So consider the following procedure:
\begin{itemize}
\item run the first SLE (i.e. $\SLE_\kappa(b)$), started from $(x,y,\dots, z_i,\dots)$ until it
reaches capacity $2\eps$.
\item  then run independently the second SLE
(i.e. $\SLE_{\tilde\kappa}(\tilde b)$) in $g_\eps^{-1}(\H)$
until it reaches capacity $2c\eps$; this capacity is measured {\em in the
original half-plane}. Let $\tilde g_{\tilde\eps}$ be the corresponding
conformal equivalence.
\item one gets two hulls resp. at $x$ and $y$ with capacity $2\eps$
and $2c\eps$; let $\phi=\tilde g_{\tilde \eps}\circ g_\eps$ be the
normalized map removing these two hulls.
\item expand $\E(F(\tilde g_{\tilde\eps}(X_\eps),\tilde Y_{\tilde \eps}))$ up to order two in $\eps$.
\end{itemize}
This describes (in distribution) how to get from $K_{0,0}$ to
$K_{\sigma,0}$, and then from $K_{\sigma,0}$ to $K_{\sigma,\tau}$. 

From the Loewner equation, it appears that $\partial_t
g'_t(w)=-2g'_t(w)/(g_t(w)-W_t)^2$. Hence
$g'_\eps(y)=1-2\eps/(y-x)^2+o(\eps)$. From the scaling property of
half-plane capacity, we get:
$$\tilde\eps=c\eps\left(1-\frac {4\eps}{(y-x)^2}\right)+o(\eps^2)$$
i.e $\tilde\eps$ is deterministic up to order two in $\eps$.
Denote by ${\mc L}$ and ${\mc M}$ the infinitesimal generators of the two SLEs:
\begin{align*}
{\mc L}&=\frac\kappa 2\partial_{xx}+b(x,y,\dots)\partial_x+\frac 2{y-x}\partial_y+\sum_{i=1}^n\frac 2{z_i-x}\partial_i\\
{\mc M}&=\frac\kappa 2\partial_{yy}+\tilde b(x,y,\dots)\partial_y+\frac 2{x-y}\partial_x+\sum_{i=1}^n\frac 2{z_i-y}\partial_i
\end{align*}
where $\partial_i=\partial_{z_i}$. Let $w=(x,y,\dots, z_i,\dots)$,
$w'=(X_\eps,g_\eps(y),\dots g_\eps(z_i),\dots )$, $w''=(\tilde
g_{\tilde\eps}(X_\eps),\tilde Y_{\tilde\eps},\dots \tilde g_{\tilde\eps}\circ
g_\eps(z_i),\dots)$. Now:
\begin{align*}
\E(F(w'')|w)&=\E(F(w'')|w'|w)=P_\eps\E(Q_{\tilde\eps}F|w')(w)\\
&=P_\eps\E\left((1+\tilde\eps {\mc M}+\frac{{\tilde\eps}^2}2 {\mc M}^2)F(w')\right)(w)
=P_\eps Q_{c\eps(1-4\eps/(y-x)^2)}F(w)+o(\eps^2)\\
&=\left(1+\eps{\mc L}+\frac{\eps^2}2{\mc
L}^2\right)\left(1+c\eps\left(1-\frac{4\eps}{(y-x)^2}\right){\mc
M}+\frac{c^2\eps^2}2{\mc M}^2\right)F(w)+o(\eps^2)\\
&=\left(1+\eps({\mc L}+c{\mc M})+\eps^2\left(\frac 12{\mc L}^2+\frac
{c^2}2{\mc M}^2+ c{\mc L}{\mc M}-\frac {4c}{(y-x)^2}{\mc
M}\right)\right)F(w)+o(\eps^2) 
\end{align*}
If we first grow a hull at $z$, then at $x$, one gets instead:
$$\left(1+\eps({\mc L}+c{\mc M})+\eps^2\left(\frac 12{\mc L}^2+\frac
{c^2}2{\mc M}^2+ c{\mc M}{\mc L}-\frac {4c}{(z-x)^2}{\mc
L}\right)\right)F(w)+o(\eps^2)$$
Hence the commutation condition reads:
\begin{equation*}
\left[{\mc L},{\mc M}\right]=\frac{4}{(y-x)^2}\left({\mc M}-{\mc L}\right)
\end{equation*}.
After simplifications, one gets:
\begin{align*}
\left[{\mc L},{\mc M}\right]+\frac{4}{(y-x)^2}\left({\mc L}-{\mc M}\right)=&
(\kappa\partial_x\tilde b-\tilde\kappa\partial_y b)\partial_{xy}\\
&+\left[\frac{2\partial_xb}{y-x}+\sum_i\frac{2\partial_ib}{y-z_i}-\tilde
b\partial_yb+\frac{2b}{(y-x)^2}+\frac{2\kappa-12}{(x-y)^3}-\frac{\tilde\kappa}2\partial_{yy}b\right]\partial_x\\
&-\left[\frac{2\partial_y\tilde b}{x-y}+\sum_i\frac{2\partial_i\tilde
b}{x-z_i}-b\partial_x\tilde b+\frac{2\tilde
b}{(x-y)^2}+\frac{2\tilde\kappa-12}{(y-x)^3}-\frac{\kappa}2\partial_{xx}\tilde
b\right]\partial_y
\end{align*}
So the commutation condition reduces to three differential
conditions involving $b$ and $\tilde b$; note the non-linear terms
$\tilde b\partial_yb $ and $b\partial_x\tilde b$.

\subsection{Rational solutions}
 
{\bf Case $n=0$:}

If $n=0$, then $b(x,y)=\rho/(x-y)$ and $\tilde
b(x,y)=\tilde\rho/(y-x)$. Then the commutation condition reduces to:
$$\left\{\begin{array}{l}
\kappa\tilde\rho=\tilde\kappa\rho\\
(\tilde\kappa-4)\rho-\rho\tilde\rho+12-2\kappa=0\\
(\kappa-4)\tilde\rho-\rho\tilde\rho+12-2\tilde\kappa=0
\end{array}\right.$$
We are only interested in the case $\kappa,\tilde\kappa>0$. Then:
$$\left\{\begin{array}{l}
\tilde\rho=\tilde\kappa\rho/\kappa\\
(\tilde\kappa-4)\rho-\rho^2\tilde\kappa/\kappa+12-2\kappa=0\\
(\kappa-4)\rho\tilde\kappa/\kappa-\rho^2\tilde\kappa/\kappa+12-2\tilde\kappa=0
\end{array}\right.$$
The last two are polynomials in $\rho$ that have a common root if and only if their resultant vanishes. This resultant (a polynomial in the coefficients) equals:
$$\frac{12\tilde\kappa(\tilde\kappa-\kappa)^2(\kappa\tilde\kappa-16)}{\kappa^3}$$
So either $\tilde\kappa=\kappa$, and then
$\rho=\tilde\rho\in\{2,\kappa-6\}$, or $\tilde\kappa=16/\kappa$, and
then $\rho=-\kappa/2$, $\tilde\rho=-\tilde\kappa/2$.

Let us comment briefly on these solutions. The condition
$\kappa\tilde\kappa=16$ obviously points at duality. The case
$\tilde\kappa=\kappa$, $\tilde\rho=\rho=\kappa-6$ corresponds in fact
to reversibility. Indeed, one has:

\begin{Lem}
An $\SLE_\kappa(\rho_1,\dots,\rho_n,\rho)$ in $\H$ started from
$(x,z_1,\dots,z_n,y)$ and aiming at $\infty$ is identical in law to a
(time-changed)
$\SLE_\kappa(\rho_1,\dots,\rho_n,\kappa-6-\rho-\sum_i\rho_i)$ in $\H$
started from $(x,z_1,\dots,z_n,\infty)$ and aiming at $y$, up to
disconnection of $y$.
\end{Lem}  
\begin{proof}
Let $(g_t)$ be the family of conformal equivalences for the first SLE,
$(K_t)$ the corresponding hulls, $W$ its driving process,
$Y^{(i)}_t=g_t(z_i)$, $Y_t=g_t(y)$. Consider the homographies:
$$\varphi_t(u)=\frac{g'_t(y)}{Y_t-u}-\frac{g''_t(y)}{2g'_t(y)}$$
and $\tilde g_t=\phi_{\varphi_0(K_t)}$. Then $\tilde
g_t\circ\varphi_0=\varphi_t\circ g_t$, and $(\tilde g_t)$ defines a
time-changed Loewner chain. Let $\tilde W_t=\varphi_t(W_t)$, $\tilde
Y^{(i)}_t=\varphi_t(\tilde Y^{(i)}_t)$, $V_t=\varphi_t(\infty)$. Then,
from It\^o's formula:
\begin{align*}
d\tilde
W_t&=d\left(\frac{g'_t(y)}{Y_t-W_t}-\frac{g''_t(y)}{2g'_t(y)}\right)=\frac{g'_t(y)}{(Y_t-W_t)}\left(-\frac
{2dt}{(Y_t-W_t)^2}-\frac{d(Y_t-W_t)}{Y_t-W_t}+\frac{\kappa
dt}{(Y_t-W_t)^2}-\frac{2dt}{(Y_t-W_t)^2}\right)\\ 
&=\frac{g'_t(y)}{(Y_t-W_t)}\left(\frac{dW_t}{Y_t-W_t}+\frac{(\kappa-6)dt}{(Y_t-W_t)^2}\right)\\
&=\frac{g'_t(y)\sqrt\kappa}{(Y_t-W_t)^2}dB_t+\frac{g'_t(y)^2}{(Y_t-W_t)^4}\left(\sum_i\rho_i\frac{Y_t-W_t}{Y^{(i)}_t-W_t}.\frac
1{\tilde W_t-V_t}+\frac{\kappa-6-\rho}{\tilde W_t-V_t}\right)dt\\ 
\end{align*}
Note that the cross-ratio is conformally invariant, so
$$\frac{Y_t-W_t}{Y^{(i)}_t-W_t}=[\infty,Y^{(i)}_t,Y_t,W_t]=[V_t,\tilde Y^{(i)}_t,\infty,\tilde W_t]=\frac{\tilde Y^{(i)}_t-V_t}{\tilde Y^{(i)}_t-\tilde W_t}$$
After a time change $ds=g'_t(y)^2dt/(Y_t-W_t)^4$, one gets:
\begin{align*}
d\tilde W_s&=\sqrt\kappa d\tilde B_s+\sum_i\rho_i\frac{\tilde
Y^{(i)}_t-V_t}{(\tilde Y^{(i)}_t-\tilde W_t)(\tilde
W_t-V_t)}dt+\frac{\kappa-6-\rho}{\tilde W_t-V_t}dt\\ 
&=\sqrt\kappa d\tilde B_s+\sum_i\frac{\rho_i}{\tilde Y^{(i)}_t-\tilde
W_t}dt+\frac{\kappa-6-\rho-\sum_i\rho_i}{\tilde W_t-V_t}dt 
\end{align*}
which defines a $\SLE_\kappa(\rho_1,\dots,\rho_n,\kappa-6-\rho-\sum_i\rho_i)$.
\end{proof}

In particular, an $\SLE_\kappa(\kappa-6)$ started from $(x,y)$ and
stopped at $\tau_y$ is simply an $\SLE_\kappa$ from $x$ to $y$. For
$\kappa=6$ and $n=0$, this is locality. 

{\bf Parametric case:}

Assume the following forms for the drift terms $b$, $\tilde b$:
\begin{align*}
b(x,y,z_i)&=\frac\rho{x-y}+\sum_i\frac{\rho_i}{x-z_i}\\
\tilde b(x,y,z_i)&=\frac{\tilde\rho}{y-x}+\sum_i\frac{\tilde\rho_i}{y-z_i}
\end{align*}
Then the commutation conditions are:
$$\left\{\begin{array}{l}
\kappa\tilde\rho=\tilde\kappa\rho\\
(\tilde\kappa-4)\rho-\rho\tilde\rho+12-2\kappa=0\\
(\kappa-4)\tilde\rho-\rho\tilde\rho+12-2\tilde\kappa=0\\
\rho\tilde\rho_i=2\rho_i\\
\tilde\rho\rho_i=2\tilde\rho_i
\end{array}\right.$$

As we have seen, the first three conditions imply that
$\kappa=\tilde\kappa$, $\rho=\tilde\rho\in\{2,\kappa-6\}$, or
$\kappa\tilde\kappa=16$, $\rho=-\kappa/2$,
$\tilde\rho=-\tilde\kappa/2$. Now, if the $\rho_i$, $\tilde\rho_i$ are
not all zero, $\rho\tilde\rho=4$, which happens if
$\rho=\tilde\rho=2$, or in the case $\kappa\tilde\kappa=16$. To sum up,
the solutions are:\\  
(i) $\kappa=\tilde\kappa$, $\rho=\tilde\rho=\kappa-6$,
$\rho_i=\tilde\rho_i=0$ (1 free parameter)\\ 
(ii) $\kappa=\tilde\kappa$, $\rho=\tilde\rho=2$, $\rho_i=\tilde\rho_i$
($(1+n)$ free parameters)\\ 
(iii) $\kappa\tilde\kappa=16$, $\rho=-\kappa/2$,
$\tilde\rho=-\tilde\kappa/2$,
$\tilde\rho_i=-(\tilde\kappa/4)\rho_i=-(4/\kappa)\rho_i$ ($(1+n)$ free
parameters)\\ 

%we study the case (ii), $\rho_i=\tilde\rho_i=2$.

These examples are ``rational'' (i.e. the drift terms are rational
functions). Yet, important examples (deduced from locality and
restriction) are transcendental. In the next section, we recast these
commutation conditions as integrability conditions, satisfied by all
these examples.

\section{Integrability conditions}

\subsection{Integrability for SLE commutation relations}

In the previous paragraph, we derived the following commutation
conditions:
$$\left\{\begin{array}{l}
\kappa\partial_x\tilde b-\tilde\kappa\partial_y b=0\\
\displaystyle\frac{2\partial_xb}{y-x}+\sum_i\frac{2\partial_ib}{y-z_i}-\tilde
b\partial_yb+\frac{2b}{(y-x)^2}+\frac{2\kappa-12}{(x-y)^3}-\frac{\tilde\kappa}2\partial_{yy}b=0\\
\displaystyle\frac{2\partial_y\tilde b}{x-y}+\sum_i\frac{2\partial_i\tilde
b}{x-z_i}-b\partial_x\tilde b+\frac{2\tilde
b}{(x-y)^2}+\frac{2\tilde\kappa-12}{(y-x)^3}-\frac{\kappa}2\partial_{xx}\tilde
b=0
\end{array}\right.$$

Now, from the first equation, one can write: 
$$b=\kappa\frac{\partial_x\psi}{\psi}, \tilde b=\tilde\kappa\frac{\partial_y\psi}{\psi}$$
for some non-vanishing function $\psi$ (at least locally). 

It turns
out that the second condition now writes:
$$-\kappa\partial_x\left(\frac{\frac{\tilde\kappa}
2\partial_{yy}\psi+\sum_i\frac{2\partial_i\psi}{z_i-y}+\frac{2\partial_x\psi}{x-y}+(1-\frac
6\kappa)\frac\psi{(x-y)^2}}{\psi}\right)=0.$$
Symmetrically, the last equation is:
$$-\tilde\kappa\partial_y\left(\frac{\frac{\kappa}
2\partial_{xx}\psi+\sum_i\frac{2\partial_i\psi}{z_i-x}+\frac{2\partial_y\psi}{y-x}+(1-\frac
6{\tilde\kappa})\frac\psi{(y-x)^2}}{\psi}\right)=0.$$
This means that a non-vanishing solution of 
$$\left\{\begin{array}{l}
\displaystyle\frac{\kappa}
2\partial_{xx}\psi+\sum_i\frac{2\partial_i\psi}{z_i-x}+\frac{2\partial_y\psi}{y-x}+\left(\left(1-\frac
6{\tilde\kappa}\right)\frac 1{(y-x)^2}+h_1(x,z)\right)\psi=0\\
\displaystyle\frac{\tilde\kappa}2\partial_{yy}\psi+\sum_i\frac{2\partial_i\psi}{z_i-y}+\frac{2\partial_x\psi}{x-y}+\left(\left(1-\frac
6\kappa\right)\frac 1{(x-y)^2}+h_2(y,z)\right)\psi=0
\end{array}\right.$$
yields drift terms $b,\tilde b$ that satisfy the commutation
condition. Obviously, these differential operators are infinitesimal
generators of the SLEs, with an added coefficient before the constant
term.

The problem is now to find functions $h_1$, $h_2$ such that
the above system has solutions (integrability conditions). 
Note that
we have not considered yet the conditions: $b,\tilde b$ translation
invariant and homogeneous of degree $(-1)$. 
This implies that $\psi$ can be chosen to be translation invariant and
homogeneous of some fixed degree.
So assume that we are given $h_1,h_2$, and a
non-vanishing (translation-invariant, homogeneous) solution $\psi$ of this system. Let:
$${\mc M}_1=\frac{\kappa}
2\partial_{xx}\psi+\sum_i\frac{2\partial_i\psi}{z_i-x}+\frac{2\partial_y\psi}{y-x}+\left(1-\frac
6{\tilde\kappa}\right)\frac 1{(y-x)^2},{\mc M}_2=\frac{\tilde\kappa}2\partial_{yy}\psi+\sum_i\frac{2\partial_i\psi}{z_i-y}+\frac{2\partial_x\psi}{x-y}+\left(1-\frac
6\kappa\right)\frac 1{(x-y)^2}$$
Then $\psi$ is annihilated by all operators in the left ideal generated by
 $({\mc M}_1+h_1)$,$({\mc M}_2+h_2)$, including in particular:
\begin{align*}
{\mc M}&=[{\mc M}_1+h_1,{\mc M}_2+h_2]+\frac 4{(x-y)^2}\left(({\mc
M}_1+h_1)-({\mc M}_2+h_2)\right)\\
&=[{\mc M}_1,{\mc M}_2]+\frac{4}{(x-y)^2}({\mc M}_1-{\mc M}_2)+([{\mc M}_1,h_2]-[{\mc M}_2,h_1])+\frac
{4(h_1-h_2)}{(x-y)^2}\\
&=\frac{3(\kappa\tilde\kappa-16)(\kappa-\tilde\kappa)}{\kappa\tilde\kappa(x-y)^4}+\left(\left(\frac{2\partial_y}{y-x}+\sum_i\frac{2\partial_i}{z_i-x}\right)h_2-\left(\frac{2\partial_x}{x-y}+\sum_i\frac{2\partial_i}{z_i-y}\right)h_1\right)+\frac
{4(h_1-h_2)}{(x-y)^2}
\end{align*}
This is an operator of order 0, so it must vanish
identically. Considering the pole at $x=y$, this implies in particular
$\tilde\kappa\in\{\kappa,16/\kappa\}$, since the fourth-order pole
must vanish. Then the second-order pole must also vanish, so
$h_1(x,z)=h(x,z)$, $h_2(y,z)=h(y,z)$ for some $h$. 
So this condition boils down to a functional equation on $h$.

For illustration, consider the following variation on an earlier
example: 
a chordal $\SLE_{8/3}$ from $x$ to $y$ is conditioned not to intersect an
independent restriction measure from $z$ to $\infty$ with index
$\nu$. Let $\varphi(x,y,z)$ be the probability of
non-intersection. Then $\varphi$ is annihilated by the operator:
$$\frac\kappa 2\partial_{xx}+\frac{\kappa-6}{x-y}\partial_x+\frac
2{y-x}\partial_y+\frac 2{z-x}\partial_z-\frac{2\nu}{(x-z)^2}$$
where $\kappa=8/3$. Obviously $\varphi$ can be expressed in terms of a
hypergeometric function. If $\psi=(y-x)^{-2\alpha}\varphi$,
$\alpha=\alpha_\kappa=5/8$, then $\psi$ is annihilated by the conjugate operators:
$$\frac\kappa 2\partial_{xx}+\frac
2{y-x}\partial_y+\frac
2{z-x}\partial_z-\frac{2\alpha}{(x-y)^2}-\frac{2\nu}{(x-z)^2},
\frac\kappa 2\partial_{yy}+\frac
2{x-y}\partial_x+\frac 2{z-y}\partial_z-\frac{2\alpha}{(y-x)^2}-\frac{2\nu}{(y-z)^2}
$$
where we also use reversibility for $\SLE_{8/3}$. It is easy to check
that in general $\tilde\kappa\in\{\kappa,16/\kappa\}$, $h(x,z)=-2\nu/(x-z)^2$
is a solution of the integrability condition
above. More generally, if $n$ points $z_1,\dots,z_n$ are marked on the
real line, a (particular) solution of the integrability condition is given by
$\tilde\kappa\in\{\kappa,16/\kappa\}$,
$$h(x,z)=-2\sum_i\frac{\mu_i}{(z_i-x)^2}-2\sum_{i<j}\nu_{ij}\left(\frac{1}{z_i-x}-\frac{1}{z_j-x}\right)^2=\sum_i\frac{\mu'_i}{(z_i-x)^2}+\sum_{i<j}\frac{\nu'_{ij}}{(z_i-x)(z_j-x)}$$
where $\mu_i,\nu_{ij}$ are real parameters. When
$\kappa=\tilde\kappa=8/3$, $\mu_i,\nu_{ij}\geq 0$, and $x<y<z_1<\cdots
z_n$, it is easy to think of a probabilistic situation corresponding
to this. Consider a chordal $\SLE_{8/3}$ from $x$ to $y$, and
condition it not to intersect independent one-sided restriction
samples $z_i\leftrightarrow\infty$ (with index $\mu_i$) and
$z_i\leftrightarrow z_j$ (with index $\nu_{i,j}$). Then reversibility
for the conditional $\SLE$ corresponds to a partition function
$\psi$ solving PDEs where $h$ is as above.

Let us get back to the functional equation for $h$:
\begin{equation}\label{inth}
\left(\left(\frac{\partial_y}{y-x}+\sum_i\frac{\partial_i}{z_i-x}\right)h(y,z)-\left(\frac{\partial_x}{x-y}+\sum_i\frac{\partial_i}{z_i-y}\right)h(x,z)\right)+2\frac
{h(x,z)-h(y,z)}{(x-y)^2}=0.
\end{equation}

We want to prove that the only solutions to this functional equation
(translation invariant and homogeneous of degree $(-2)$) are the
rational functions given above when there are at most $3$ marked $z$
points (including infinity).  
By expanding in $\eps$ where $y=x+\eps$, one sees that $h$ must be
annihilated by the family of operators:
$$\ell_{0,n}=\frac n{(n+1)(n+2)}\partial_x^{n+2}+\sum_i\left(\frac{\partial_i\partial_x^n}{z_i-x}-\frac{n!\partial_i}{(z_i-x)^{n+1}}\right)$$
for $n\geq 0$. Also, $h$ must be translation invariant and homogeneous
of degree $-2$. So for $n=2$, one can write $h(x,y,z)=\tilde
h((y-x)/(z-x))/(z-x)^2$ for instance, and $\tilde h$ satisfies a
third-order ODE (since $\ell_{0,1}h=0$); but we have already exhibited 3 linearly independent
solutions in this case, so the classification is complete for $n=1,2$.
In particular, in the case $n=2$, this is closely related to the
discussion in Section 8.5 of \cite{LSW3}.

When $n\geq 3$, the configuration $(\H,z_1,\dots,z_n,\infty)$
corresponds to a $(n-2)$-dimensional moduli space. We already know $n(n+1)/2$
(linearly independent) rational solutions. We proceed to show that
arbitrary (smooth) functions on this ``residual moduli space'' lead to
solutions of the functional equation.

Define $\ell_x$, $\ell_y$ to be the differential operators:
$$\ell_x=\frac 2{y-x}\partial_y+\sum_i\frac 2{z_i-x}\partial_i,{\rm\ \
\ }\ell_y=\frac 2{x-y}\partial_x+\sum_i\frac 2{z_i-y}\partial_i$$
representing the Loewner flow with singularities at $x,y$
respectively, restricted to marked points. Now the operators
$$\tilde\ell_x=\frac 6{y-x}\partial_x+\ell_x,{\rm\ \
\ }\tilde\ell_y=\frac 6{x-y}\partial_y+\ell_y$$
are the generators of $\SLE_0(-6)$, i.e. the hyperbolic geodesic from
$x$ to $y$. In this case, we have seen that the commutation relation:
$$[\tilde\ell_x,\tilde\ell_y]=\frac 4{(x-y)^2}(\tilde\ell_y-\tilde\ell_x)$$
is satisfied. Hence if $f$ is a positive, translation invariant, homogeneous of
degree 0 function of $z_1,\dots,z_n$ (i.e. a function on the
``residual moduli space''), then the weight:
$$h(x,z)=\tilde \ell_x\log f(z)=\ell_x\log f(z)$$
satisfies the functional equation (\ref{inth}) and is translation
invariant and homogeneous of degree $(-2)$.

Prompted by these solutions, we want to get a global interpretation of
the differential condition (\ref{inth}). Consider two simple paths
$\gamma$, $\tilde\gamma$ in $\H$ started from $x,y$ respectively, non
intersecting and driven by smooth functions, for simplicity. The union
of the two paths can be parametrized by Loewner's equation in two
extremal fashions: exploring $\gamma$ first and then $\tilde\gamma$,
or exploring $\tilde\gamma$ first and then $\gamma$. Define:
$$c(\gamma,z)=\int_0^\tau h(x_t,z_1(t),\dots,z_n(t))dt$$
where $\gamma$ has half-plane capacity $2\tau$, and $z_1(.),\dots,z_n(.)$ follow the
Loewner flow driven by $(x_t)$. If $\gamma$, $\tilde\gamma$ have half-plane capacity
$\eps,c\eps$ respectively, then:
$$c(\gamma,z)+c(\phi_{\gamma}(\tilde\gamma),\phi_\gamma(z))=
c(\tilde\gamma,z)+c(\phi_{\tilde\gamma}(\gamma),\phi_{\tilde\gamma}(z))+O(\eps^2).$$
This is exactly the content of (\ref{inth}), the corrections being
as in the commutation relations. So it is easy to see (and
straightforward to write, using comparisons of Loewner chains etc...)
that for macroscopic paths $\gamma$, $\tilde\gamma$, we have:
$$c(\gamma,z)+c(\phi_{\gamma}(\tilde\gamma),\phi_\gamma(z))=
c(\tilde\gamma,z)+c(\phi_{\tilde\gamma}(\gamma),\phi_{\tilde\gamma}(z)).$$
Then we can define $c(\gamma\cup\tilde\gamma,z)$ to be this
quantity. Similarly, one can grow a first half of $\gamma$, then a
half of $\tilde\gamma$, then grow the end of $\gamma$,\dots, and get the same quantity. If $A$ is any smooth hull (not intersecting the $z_i$'s),
one can define $c(A,z)=c(\partial A,z)$, where $\partial A$ is a
smooth curve. Notice that this does not depend on the orientation of
$\partial A$. Indeed, let $x,y$ be the endpoints of $\partial A$,
$\gamma$, $\tilde\gamma$ as above, such that $\gamma$ and
$\tilde\gamma$ are at distance at most $\eps$ of $\partial A$ in the
Hausdorff metric. Then $\phi_\gamma(\tilde\gamma)$ and
$\phi_{\tilde\gamma}(\gamma)$ have small half-plane capacity, and one
can apply the previous result. If the boundary $\partial A$ is
described as the union of two arcs starting at $x,y$ respectively, one
gets the same result.
Also, if $A,B$ are two smooth hulls
contained in a compact hull not intersecting the $z_i$'s, then, by
similar arguments, we have
the Lipschitz condition:
$$|c(A,z)-c(B,z)|\leq k({\rm cap}(A\cup B)-{\rm cap}(A\cap B)).$$
Hence we can extend the definition of $c$ by approximation by smooth hulls.

So if $h$ is a function solving (\ref{inth}), then we can define a
function $C=\exp c$ on hulls and residual configurations such that:
\begin{enumerate}
\item For all hulls $A,B$,
 $C(A.B,z)=C(B,\phi_A(z))C(A,z)$.
\item If $A$ is a hull of half-plane capacity $\eps$ located at $x$, then
$C(A,z)=1+2\eps h(x,z)+o(\eps)$. 
\end{enumerate}

Here $A.B$ designates the concatenation of the two hulls $A,B$:
$\phi_{A.B}=\phi_B\circ\phi_A$.
Conversely, if a nice function $C$ satisfying (i) is given, we can
recover a function $h$ satisfying (\ref{inth}) (its derivative in the
direction $\ell_x$). So far, we have considered the following $C$'s:
$$C(A,z_1,\dots,z_{n+1})=\prod_{i<j}\left(\phi'_A(z_i)\phi'_A(z_j)\left(\frac{z_j-z_i}{\phi_A(z_j)-\phi_A(z_i)}\right)^2\right)^{\nu_{ij}}\frac
{f(\phi_A(z_1),\dots,\phi_A(z_{n+1}))}{f(z_1,\dots,z_{n+1})}$$
where the marked point at infinity is now $z_{n+1}$. All the factors
here are conformally invariant.

Let ${\mc H}$ designates a semigroup of hulls (as in
\cite{LSW3}). The residual configuration is
$(\H,z_1,\dots,z_{n+1}=\infty)$ (it no longer depends on the marked
points $x,y$).
More precisely, an element of ${\mc H}$ is a compact
subset $A$ of $\overline\H$ such that $\H\setminus A$ is simply-
connected and $\partial A\cap\R\subset \overline{A\cap\H}$; the
semigroup is concatenation: $\phi_{A.B}=\phi_B\circ\phi_A$.
Let ${\mc F}$ be a vector space of 
functions of the residual configuration $(z_1,\dots,z_{n+1})$ (say
smooth functions). Then ${\mc H}$ acts on ${\mc F}$ by
$A.f=f\circ\phi_A$ (note that this is not everywhere defined). If $C$
is as above, $C$ can be seen as a map ${\mc H}\rightarrow{\mc F}$ by
$C(A)(z)=C(A,z)$. Then the condition (i) reads:
$$C(A.B)=C(A)(A.C(B))$$
which is saying that $C$ is a 1-cocycle for (multiplicative) group
cohomology for the ${\mc H}$-module ${\mc F}$. This is formal since
${\mc H}$ is only a semigroup and the operation is not everywhere
defined. Then the question is to determine the first cohomology group
$H^1({\mc H},{\mc F})$, restricted to cocycles with M\"obius covariance:
$$C(\phi(A),\phi(z))=C(A,z)$$
where $\phi$ an homography. Let $c_{ij}$ be the cocycle:
$$c_{ij}(A)(z)=\phi'_A(z_i)\phi'_A(z_j)\left(\frac{z_j-z_i}{\phi_A(z_j)-\phi_A(z_i)}\right)^2.$$
We have a natural map:
\begin{align*}
\R^{n(n+1)/2}&\longrightarrow H^1({\mc H},{\mc F})\\
(\nu_{ij})_{i<j}&\longmapsto \prod_{i<j}(c_{ij})^{\nu_{ij}}
\end{align*}
The content of Section 8.5 in \cite{LSW3} is that this is an
isomorphism when $n\leq 2$. At this point, it is not quite clear
whether this is onto in general. Though it not one-to-one as soon as $n\geq
3$, since e.g. $c_{12}c_{34}/c_{13}c_{24}$ is a coboundary. It is easy
to see that its image has dimension $n+1$ and is generated for
instance by $(c_{12},c_{13},\dots,c_{1,n+1},c_{23})$.

So let us extend the discussion in \cite{LSW3} when $n\geq
3$. This involves several complications. The idea is to restrict to a
subgroup of ${\mc H}$ (or rather a subalgebra of its tangent algebra)
fixing a configuration $(z_1,\dots,z_{n+1})$. A cocyle restricts to a
character of this subgroup, that must vanish on commutators. This
gives differential conditions for tangent cocycles, that can then be integrated.

We begin by fixing three marked points at $0,1,\infty$ (which we
can do by a M\"obius transformation). The other marked points are $z_1,\dots,z_{n-2}$.
Consider the real Lie algebra generated by the Loewner fields:
$$A(x)\partial_w=\frac{w(w-1)}{w-x}\partial_w$$
that fix $0,1,\infty$.
Consider the subalgebra fixing the marked points
$z_1,\dots,z_{n-2}$. This includes the fields:
$$A(x_1,\dots,x_{n-1})=\det(A(x_i)(z_j))_{1\leq i,j\leq n-1}=\prod_{i=1}^{n-1}z_i(z_i-1)\frac{V(x_1,\dots,x_{n-1})V(z_1,\dots,z_{n-1})}{\prod_{i,j}z_i-x_j}$$
with the convention $w=z_{n-1}$, where $V$ is the Vandermonde
polynomial $V(x_1,\dots,x_{n-1})=\prod_{i<j}(x_j-x_i)$ (Cauchy determinant).
By taking limits of such vector fields (as $x_1,\dots,x_n\rightarrow
x$, normalizing by $V(x_1,\dots,x_n)$), one gets
$\tilde A(x)=(n-2)!P(w)/(w-x)^{n-1}$, where
$P(w)=w(w-1)(w-z_1)\cdots(w-z_{n-2})$ is a polynomial of
degree $n$.
The (tangent) cocycle $dC$ restricted to this subalgebra is a morphism
to the (trivial) Lie algebra $\R$, so it vanishes on commutators.
By considering the limit of $[\tilde A(x),\tilde A(y)]$ as $y\searrow
x$, we see that $\hat A(x)=(2n-1)!P(w)^2/(w-x)^{2n}$ is annihilated by
$dC$. It is easy to see that:
$$\hat A(x)=\left(\partial_x^{2n-1} P(x)\partial_x^{2-n}\right)\tilde A(x).$$
Note that $(\partial_x^{2n-1} P(x)\partial_x^{2-n})$ is a
well-defined linear differential operator of degree $n+1$, with
polynomial coefficients. It follows that:
$$dC(\tilde A(x))=a+\frac{a_0}{x^{n-1}}+\frac{a_1}{(x-1)^{n-1}}+\sum_{j=1}^{n-2}\frac{b_j}{(x-z_j)^{n-1}}$$
where the coefficients depend on the marked points
$z_1,\dots,z_{n-2}$. We want to integrate this to get information on
$dC(A(x_1,\dots))$. 
%First, we have,
%$$(\partial_{x_2}\partial_{x_3}^2\dots\partial_{x_{n-1}}^{n-2})_{|x_1=\cdots
%=x}A(x_1,\dots,x_{n-1})=1!\dots (n-2
%)!V(z_1,\dots,z_{n-2})\prod_{i=1}^{n-2}\frac{z_i(z_i-1)}{(z_i-x)^{n-1}}\tilde A(x).$$
Consider the atomic measure (in the variable $x$):
$$m=\sum_i(-1)^{n-1-i}V(x_1,\dots,\hat {x_i},\dots)\delta_{x_i}.$$
One can think of $m$ as the Vandermonde determinant $V(x_1,\dots,x_{n-1})$
where the last row $(x_1^{n-2},\dots,x_{n-1}^{n-2})$ is replaced
by $(\delta_{x_1},\dots,\delta_{x_{n-1}})$.
Then we have the decomposition:
$$\frac{V(x_1,\dots,x_{n-1})}{\prod_i (w-x_i)}=\left\langle m,\frac
1{w-x}\right\rangle=(-1)^n\left\langle m^{(2-n)},\frac
{(n-2)!}{(w-x)^{n-1}}\right\rangle$$
where $m^{(2-n)}$ can be chosen with a continuous, compactly supported
density (since $\langle m,1\rangle=\cdots=\langle m,x^{n-3}\rangle=0$,
which follows from the Vandermonde form of $m$). It follows that:
$$A(x_1,\dots,x_{n-1})=(-1)^n\prod_{i=1}^{n-2}z_i(z_i-1)\frac{V(z_1,\dots,z_{n-2})}{\prod_{i<n-1,j\leq
n-1}z_i-x_j}\left\langle m^{(2-n)},\tilde A(x)\right\rangle$$
and consequently:
$$dC(A(x_1,\dots,x_{n-1}))=\prod_{i=1}^{n-2}z_i(z_i-1)\frac{V(z_1,\dots,z_{n-2})}{\prod_{i<n-1,j\leq
n-1}z_i-x_j}\left\langle m,ax^{n-2}+\frac{a_0}{x}+\frac{a_1}{x-1}+\sum_{j=1}^{n-2}\frac{b_j}{x-z_j}\right\rangle$$
for some coefficients $a,a_0,b_j$ depending on
$z_1,\dots,z_{n-2}$. From the Vandermonde form of $m$ and simple
manipulations, we get:
$$dC(A(x_1,\dots,x_{n-1}))=\prod_{i=1}^{n-2}z_i(z_i-1)\frac{V(z_1,\dots,z_{n-2})V(x_1,\dots,x_{n-1})}{\prod_{i<n-1,j\leq
n-1}z_i-x_j}\left(a+\frac{a_0}{\prod_jx_j}+\frac{a_1}{\prod_jx_j-1}+\sum_{k=1}^{n-2}\frac{b_k}{\prod_jx_j-z_k}\right).$$
Using the Cauchy determinant formula (and its derivative w.r.t $z_k$), it is then
easy to see that:
$$\displaystyle{
\left|\begin{array}{ccc}
\frac 1{x_1-z_1}&\dots&\frac 1{x_{n-1}-z_1}\\
\vdots&&\vdots\\
\frac 1{x_1-z_{n-2}}&\dots&\frac 1{x_{n-1}-z_{n-2}}\\
R(x_1)(z)&\dots&R(x_{n-1})(z)
\end{array}\right|=0}$$
where:
$$R(x)=dC(A(x))-\left(a+\frac{a_0}{x}+\frac{a_1}{x-1}+\sum_{j=1}^{n-2}\frac{b_j}{(x-z_j)^2}\right)$$
which depends implicitly on the $z$ variables (and we have replaced
$a_0$ with $a_0/\prod_i z_i$, $a_1$ with $a_1/\prod_i (z_i-1)$; since
we use $0,1$ for normalization, this induces some asymmetry).
For the part corresponding to the coefficient $a$, one can
consider the limit $w\rightarrow\infty$ of $A(x_1,\dots,x_{n-1}$).
If $Q(x)=(x-z_1)\dots (x-z_{n-2})$, after row operations, we get:
$$
\left|\begin{array}{ccc}
1&\dots&1\\
\vdots&&\vdots\\
x_1^{n-3}&\dots& x_{n-1}^{n-3}\\
Q(x_1)R(x_1)(z)&\dots&Q(x_n)R(x_{n-1})(z)
\end{array}\right|=0$$
If we apply $\partial_{x_1}^{n-2}$ to this
expression (using multilinearity in rows), we get
$\partial_{x_1}^{n-2}Q(x_1)R(x_1)(z)=0$. Hence:
$$dC(A(x))=a+\frac{a_0}{x}+\frac{a_1}{x-1}+\sum_{j=1}^{n-2}\frac{c_j}{x-z_j}+\sum_{j=1}^{n-2}\frac{b_j}{(x-z_j)^2}$$
where coefficients may depend on $z$. 

In terms of $h$, this means that:
$$h(x,z)=\sum_i\frac{\mu_i(z)}{(x-z_i)^2}+\sum_i\frac{\rho_i(z)}{x-z_i}$$
where $\mu_i,\rho_i$ are translation invariant and homogeneous of
degree $0,-1$ respectively (for normalization reasons, we have to
divide $dC(A(x))$ by $x(x-1)$). 
So we can substitute this expression in the commutation equation (\ref{inth}), and get
after some simplifications:
$$\sum_{i,j}\left(\frac 1{(y-z_j)(x-z_i)^2}-\frac 1{(x-z_j)(y-z_i)^2}\right)\partial_j\mu_i+\left(\frac 1{(y-z_j)(x-z_i)}-\frac 1{(x-z_j)(y-z_i)}\right)\partial_j\rho_i=0$$
as a rational function in $x,y$. This implies that
$\partial_j\mu_i=0$ for all $i,j$ (considering first the coefficient
of $(x-z_i)^{-2}$, and then letting $y$ vary). Similarly, considering
the coefficient of $(x-z_i)^{-1}$, we get:
$$\sum_j\frac{\partial_j\rho_i-\partial_i\rho_i}{y-z_j}
=0$$
and letting $y$ vary, we get the cross-derivative condition
$\partial_i\rho_j=\partial_j\rho_i$.
Hence we can find constant coefficients $\mu_i$ and a function $f(z)$ such that:
$$h(x,z)=\sum_i\frac{\mu_i}{(x-z_i)^2}+\ell_x\log f(z).$$
Note that $f(z)=(z_j-z_i)^{\nu_{ij}}$ produces the term
$2\nu_{ij}/(x-z_i)(x-z_j)$. So we can rewrite $h$ as:
$$h(x,z)=\sum_{i,j}\nu_{ij}\left(\frac 1{x-z_i}-\frac 1{x-z_j}\right)^2+\ell_x\log f(z)$$
in a non-unique fashion (e.g. take $f$ the power of a cross-ratio; here $z_{n+1}=\infty$). Since:
$$\partial_j\sum_i\rho_i=(\sum_i\partial_i)\rho_j=0,{\rm\ \ }
\partial_j\sum_iz_i\rho_i=(\sum_iz_i\partial_i\rho_j)+\rho_j=0$$
$f$ is translation invariant and homogeneous of some fixed
degree. This degree can be set to 0 by adjusting the $\nu_{ij}$'s.
One can freely set all the $\nu_{ij}$'s to zero except
$(\nu_{12},\dots, \nu_{1,{n+1}},\nu_{23})$ e.g., which gives $(n+3)$
numerical invariants.

Let us sum up the previous discussion.

\begin{Lem}
If the system
$$\left\{\begin{array}{l}
\displaystyle\frac{\kappa}
2\partial_{xx}\psi+\sum_i\frac{2\partial_i\psi}{z_i-x}+\frac{2\partial_y\psi}{y-x}+\left(\left(1-\frac
6{\tilde\kappa}\right)\frac 1{(y-x)^2}+h_1(x,z)\right)\psi=0\\
\displaystyle\frac{\tilde\kappa}2\partial_{yy}\psi+\sum_i\frac{2\partial_i\psi}{z_i-y}+\frac{2\partial_x\psi}{x-y}+\left(\left(1-\frac
6\kappa\right)\frac 1{(x-y)^2}+h_2(y,z)\right)\psi=0
\end{array}\right.$$
admits a non-vanishing solution $\psi$ (smooth, homogeneous and translation
invariant), then:
\begin{enumerate}
\item Either $\tilde\kappa=\kappa$ or $\tilde\kappa=16/\kappa$.
\item The functions $h_1$, $h_2$ can be written as $h_1(x,z)=h(x,z)$,
$h_2(y,z)=h(y,z)$, where:
$$h(x,z)=\sum_{1\leq i<j\leq n+1}2\nu_{ij}\left(\frac 1{x-z_i}-\frac 1{x-z_j}\right)^2+\ell_x\log f(z)$$
where $\nu_{ij}$ are constant parameters and $f(z)$ is a conformally
invariant function of the marked points $(z_1,\dots,z_{n+1})$.  
\end{enumerate}
Moreover the weight vector $(\nu_{ij})_{i<j}$ is well-defined modulo
the relations $\nu_{ij}+\nu_{kl}=\nu_{il}+\nu_{jk}$.
\end{Lem}

Conversely, it is natural to ask whether these conditions are
sufficient for the existence of a non vanishing solution $\psi$.
In the situation where $h$ is the rational function:
$$h(x,z)=\sum_i\frac{\mu_i}{(x-z_i)^2}+\sum_{i<j}\frac{\nu_{ij}}{(x-z_i)(x-z_j)}$$
then we can find an elementary solution
of the form:
$$\psi(x,y,z)=\left(\prod_i(z_i-x)^{a_i}(z_i-y)^{\tilde a_i}\right)(y-x)^b\prod_{i<j}(z_j-z_i)^{b_{ij}}$$
If $\tilde\kappa=\kappa$, we can pick $b=2/\kappa$, $\tilde a_i=a_i$,
$b_{ij}=\nu_{ij}/2$ and $a_i$ a solution of
$(\kappa/2)a_i(a_i-1)+2a_i+\mu_i=0$ for all $i$. If
$\tilde\kappa=16/\kappa$, set $b=-1/2$, $\tilde a_i=-a_i\kappa/4$, $b_{i,j}=\nu_{ij}/2$ and $a_i$ a solution of
$(\kappa/2)a_i(a_i-1)+2a_i+\mu_i=0$ for all $i$. 
Finding $\psi$ of a prescribed homogeneity degree is much more
difficult.

For a general $h$, we shall later discuss martingale interpretations
of solutions.
If one considers $n$ SLEs, $n\geq 2$, one gets a system of $n$ linear
PDEs (with coefficients $h_1,\dots h_n$ to be specified). We give an
example of this situation in the next section.

\subsection{Restriction, locality and cocycles}

In the previous subsection, we made some progress on the question of classifying
restriction measures, which we now explicit (in the chordal setting).
A configuration is a simply connected domain $D$ with $n$ marked
points on the boundary: $(D,z_1,\dots,z_n)$ (say $D$ has Jordan
boundary). The following definition is a natural extension of \cite{LSW3}. We call a restriction measure a
collection $\mu$ of measures parametrized by the configuration
$(D,z_1,\dots,z_n)$ such that:
\begin{enumerate}
\item The measure $\mu_{(D,z_1,\dots)}$ is supported on simply-
connected compact
subsets of the compactification of $D$ that intersect the boundary of
$D$ exactly at $z_1,\dots,z_n$. Also, $\partial K\cap\partial D\subset \overline{K\cap D}$.
\item (Conformal invariance) If $\phi:(D,z_1,\dots)\rightarrow(D',z_1',\dots)$ is an
equivalence of configurations, $\phi_*\mu_{(D,\dots)}=\mu_{(D',\dots)}$.
\item (Restriction property) If $(D,z_1,\dots)$ and $(D',z_1,\dots)$ are configurations such
that $D'\subset D$ and $\partial D'\cap\partial D$ contains
neighbourhoods of $z_1,\dots,z_n$, then:
$\mu_D(.|.\subset D')=\mu_{D'}$.
\end{enumerate}

Given such a restriction measure (strictly speaking, collection of
measures, or measure-valued function on the moduli space), we can
define a cocycle: $C_\mu(A,z)=\mu_{(\H,z)}(\{K: K\cap A=\varnothing\})$.
Indeed, the tower property of conditional expectations:
$$\mu_{(\H,z)}(.|\subset \H\setminus A.B)=
\mu_{(\H,z)}(.|\subset \H\setminus A.B|\subset \H\setminus A)
=\mu_{(\H\setminus A,z)}(.|\subset \H\setminus A.B)
=\mu_{(\H,A.z)}(.|\subset \H\setminus B)$$
translates into the cocycle condition: $C_\mu(A.B,z)=C_\mu(A,z)C_\mu(B,A.z)$.
Besides, $\mu$ is entirely determined by $C_\mu$, since the events
$\{K\cap A=\varnothing\}$ constitute a $\pi$-system generating the full
$\sigma$-algebra for random compact sets $K$ satisfying the
topological condition (i).

We have proved that under a regularity assumption $C_\mu$ can be
expressed as:
$$C_\mu(A,z)=\prod_{i<j}\left(\phi'_A(z_i)\phi'_A(z_j)\left(\frac{z_j-z_i}{\phi_A(z_j)-\phi_A(z_i)}\right)^2\right)^{\nu_{ij}}\frac
{f(\phi_A(z_1),\dots,\phi_A(z_n))}{f(z_1,\dots,z_n)}.$$
The question is what cocycles can be realized through a restriction
measure. A result of \cite{LSW3} is that if $n=2$, then the cocycle
$\left(\phi'_A(z_1)\phi'_A(z_2)\left(\frac{z_2-z_1}{\phi_A(z_2)-\phi_A(z_1)}\right)^2\right)^{\nu}$
corresponds to a restriction measure iff $\nu\geq 5/8$. Also, there is
an operation $\bullet$ on restriction measures (filling of the union of
independent samples) such that $C_{\mu_1\bullet\mu_2}=C_{\mu_1}C_{\mu_2}$.

Let us now discuss locality; the close relation between locality and
restriction is stressed in \cite{LSW3}.
Consider a configuration $(D,a,b,c)$ with three marked points on the
boundary, in cyclic order.
Using chordal $\SLE_6$, one can define a distribution $\mu$ relative
to such configurations and satisfying:
\begin{enumerate}
\item $\mu_{(D,a,b,c)}$ is supported on simply connected compact subsets $K$ of the
compactification of $D$ whose intersection with $\partial D$ consists
in a point $X\in (bc)$ and an arc contained in $(cb)$ and containing
$a$. Also, $\partial K\cap\partial D\subset \overline{K\cap D}$.
\item (Conformal invariance) If $\phi:(D,a,b,c)\rightarrow
(D',a',b',c')$ is an equivalence of configurations, then 
$\phi_*\mu_{(D,a,b,c)}=\mu_{(D',a',b',c')}$.
\item (Locality) Let $\mu_{(D,a,b,c,x)}$ be the disintegrated measure 
$\mu_{(D,a,b,c)}(.|X\in dx)$. 
If $(D,a,b,c)$ and $(D',a,b',c')$ are configurations such
that $D'\subset D$, $b'\in (ab)$, $c'\in(ca)$, $\partial
D'\cap\partial D$ contains the arcs $(ab')$, $(c'a)$ and
a neighbourhood of $x$, then:
$$\mu_{(D,a,b,c,x)}(.|K\subset D',\partial K\subset (ab')\cup\{x\}\cup
(c'a))=\mu_{(D',a,b',c',x)}.$$
\end{enumerate}
Let us associate a cocycle to this collection of measures. We phrased
the locality property as a restriction property for disintegrated
measures, so we can define as above:
$$C(A,z)=\mu_{(\H,a,b,c,x)}(\{K:K\subset \H\setminus A\})$$
where $z=(a,b,c,x)$. With these conditions, the measures are entirely
determined by the distribution of $X$. More precisely, if this
distribution is $f(x)dx$ for the configuration $(D,a,b,c)$ and $\phi$
is the equivalence $(D,a,b,c)\rightarrow (D',a,b',c')$, then
$C_\mu(A,z)=\phi'(x) f\circ\phi(x)/f(x)$, or, without assumptions on the
normalization of the equivalence $\phi$:
$$C_\mu(\phi,z)=\frac{c_{ax}c_{bc}^{1/2}}{c_{ab}^{1/2}c_{ac}^{1/2}}(\phi,z)\frac {f\circ\phi(z)}{f(z)}$$
and again one can recover $\mu$ from $C_\mu$. One can think of several
extensions to configurations with more points. Consider the following
example, coming from percolation. For critical percolation in a
rectangle $(ABCD)$, consider the highest (resp. lowest) open path from
$(AB)$ to $(CD)$ with endpoints $A',B',C',D'$ respectively. This gives
a measure (with mass $<1$) on hulls delimited by the two
paths. Disintegrated measures (w.r.t. $A',\dots,D'$) have the
restriction property. The measure on hulls is determined by the joint
distribution $(A',B',C',D')$, than can be obtained by taking partial
derivatives of Cardy's formula (see \cite{D5} for related questions).

\section{A particular case}

In this section, we discuss the important situation where all marked
points on the boundary are growth points for commuting SLEs. This 
situation is studied in greater details in \cite{D7}.

As described earlier, consider the half-plane $\H$ with
$(2n)$ marked points on the boundary, $(x_1,\dots, x_{2n})$, in
cyclical order. Consider $n$ independent
$\SLE_{8/3}$ from $x_{2i-1}$ to $x_{2i}$, $i=1,\dots,n$. Define
$\psi(x_1,\dots,x_{2n})$ to be the probability of no pairwise
intersection. This function is
invariant under the full M\"obius group ($\infty$ is only used for
normalization). Now, if $\gamma_i$ is the trace of the $i$-th SLE, define 
\begin{align*}
\psi(x_i)
&=\P(\gamma_i\cap\gamma_j=\varnothing, 1\leq i< j\leq n)\\
&=\P(\gamma_1(0,t)\cap(\cup_{j>1}\gamma_j)=\varnothing,\gamma_1(t,\infty)\cap(\cup_{j>1}\gamma_j)=\varnothing,
\gamma_i\cap\gamma_j=\varnothing, 2\leq i< j\leq n)\\
&=\prod_{j>1}\P(\gamma_j\cap\gamma_1(0,t)=\varnothing)\P(\tilde\gamma_i\cap\tilde\gamma_j=\varnothing, 1\leq i< j\leq n)
\end{align*}
where $\tilde\gamma_i$ are independent $\SLE_{8/3}$ in the domain
$(\H,g_t(\gamma_1(t)),g_t(x_2),\dots,g_t(x_{2n}))$ and $(g_t)$ are the
conformal equivalences associated with the first SLE. This relies on
the Markov property for $\gamma_1$, the restriction property for each
$\gamma_j$, $j>1$, and induction on $n$. As a consequence, the following process:
$$\psi(W_t,g_t(x_2),\dots,g_t(x_{2n}))\prod_{j>1}\left(g'_t(x_{2j-1})g'_t(x_{2j})\left(\frac{x_{2j}-x_{2j-1}}{g_t(x_{2j})-g_t(x_{2j-1})}\right)^2\right)^{\alpha_\kappa}$$
is a martingale, where $\alpha_\kappa=(6-\kappa)/2\kappa$, $\kappa=8/3$. Now, one
can do this starting at each point $x_i$ (since $\SLE_{8/3}$ is
reversible). This implies that $\psi$ is annihilated by the operators
($k=1,\dots, 2n$):
$$\frac\kappa 2\partial_{kk}+\sum_{l\neq k}\frac{2\partial_{l}}{x_l-x_k}+(\kappa-6)\frac{\partial_{k}}{x_{k}-x_{\iota(k)}}+\frac{\kappa-6}\kappa\sum_{\{j,\iota(j)\}\neq
\{k,\iota(k)\}}\left(\frac1{x_{j}-x_{k}}-\frac1{x_{\iota(j)}-x_{k}}\right)^2
$$
where $\iota$ defines the chosen pairing $\iota(2i-1)=2i$, $\iota(2i)=2i-1$.

This is not very symmetrical. It is easy to see that the function
$$\psi(\dots x_i\dots)\prod_j(x_{2j}-x_{2j-1})^{1-6/\kappa}$$
is annihilated by the operators:
$$\left\{\begin{array}{l}\displaystyle\frac\kappa 2\partial_{k,k}+\sum_{l\neq
k}\frac{2\partial_{l}}{x_l-x_k}+\frac{\kappa-6}\kappa\sum_{l\neq
k}\frac 1{(x_l-x_k)^2}, \hspace{1cm} k=1,\dots, 2n\\
\sum_k\partial_k\\
\sum_k x_k\partial_k-n(1-6/\kappa)\\
\sum_k x_k^2\partial_k-(1-6/\kappa)(x_1+\cdots+ x_{2n})
\end{array}\right.$$
the last three ones corresponding to the invariance of $\psi$ under the Moebius group.
In fact, as is discussed in \cite{D7}, one can make sense of this
sytem for any $\kappa\in (0,8/3)$ using appropriate loop-soups.

Let us make a few remarks on this system. First, each choice of a
non-crossing pairing of the $(2n)$ boundary points yields a solution;
there are $C_n$ such pairings. If $\kappa=6$, this is the system
satisfied by crossing probabilities for critical percolation in a
$(2n)$-gon with alternating boundary conditions. The number of these
crossing probabilities is the number of non-crossing partitions of the
set of blue edges, which is known to be $C_n$. In the case $n=2$, it
is trivial to solve this system, which reduces to a hypergeometric
equation (and $C_2=2$). If $n=3$, one can write this system in a
Pfaffian form, proving that its rank is indeed $C_3=5$. In the case
$\kappa=6$, $n=3$, and configurations with 3-fold symmetry, one can
express solutions in terms of $_3F_2$. Finally, one
can take the limit $\kappa\rightarrow\infty$ of the system; in this
case, solutions are polynomials, and it is easy to see that the rank
of the system is $C_n$ for all $n$. Euler integrals for solutions of
this system are discussed in \cite{D7}.

\section{Local commutation}

In this section we see how to go from infinitesimal commutation
relations to commutation (in law) of SLE hulls. Recall from Section 3
the definition of an $\SLE(\kappa,b,\tilde\kappa,\tilde b)$. We have seen in the
previous sections that the existence of such an SLE implies conditions
on $(\kappa,b,\tilde\kappa,\tilde b)$ (in particular either
$\tilde\kappa=\kappa$ or $\tilde\kappa=16/\kappa$). Conversely, assume
that the data $(\kappa,b,\tilde\kappa,\tilde b)$ satisfies the
appropriate conditions. We will see that this implies the existence of
an $\SLE(\kappa,b,\tilde\kappa,\tilde b)$. Note that this is not
saying anything on the long time behaviour of such an SLE. The
questions involving collisions of commuting SLEs are delicate
and cannot be handled by these methods.

\begin{Prop}
Consider the upper half-plane $\H$ with $n+2$ distinct marked real
points $x,y,z_1,\dots,z_n$. Assume that $\kappa,\tilde\kappa$ are
positive numbers and $b,\tilde b$ are smooth functions, translation
invariant and homogeneous of degree $(-1)$, such that the following
relation holds:
\begin{equation*}
\left[{\mc L},{\mc M}\right]=\frac{4}{(y-x)^2}\left({\mc M}-{\mc L}\right)
\end{equation*}
where ${\mc L}$ (resp. ${\mc M}$) is the infinitesimal generator of
$\SLE_\kappa(b)$ (resp. $\SLE_{\tilde\kappa}(\tilde b)$) growing at
$x$ (resp. $y$). Then there exists an
$\SLE(\kappa,b,\tilde\kappa,\tilde b)$.
\end{Prop}

We will use the following lemma. 
Let $\phi_0$ be some conformal
equivalence $\H\setminus K\rightarrow\H$, with hydrodynamic
normalization at infinity. Let us call $\phi_0$-capacity the increasing function
on hulls: ${\rm cap}\circ\phi_0^{-1}-{\rm cap}(K)$.

\begin{Lem} 
With the hypotheses of the Proposition, 
let $D_1$, $D_2$ be disjoint compact neighbourhoods of $x$, $y$ resp.,
with Jordan boundary, not containing any other marked point; $\eta_1$,
$\eta_2$ are two positive numbers.
Then the two following procedures define the same
probability law on pairs of chains in $\H$:
\begin{enumerate} 
\item Grow an $\SLE_\kappa(b)$ at $x$ until it exits $D_1$ 
or its $\phi_0$-capacity exceeds $\eta_1$, 
Then grow an independent
$\SLE_{\tilde\kappa}(\tilde b)$ at $y$ in the remaining domain,
until it exits $D_2$ or its $\phi_0$-capacity exceeds $\eta_2$.
\item Grow an $\SLE_{\tilde\kappa}(\tilde b)$ at $y$ until it exits
$D_2$ 
 or its $\phi_0$-capacity exceeds $\eta_2$. 
Then grow an independent
$\SLE_\kappa(b)$ at $x$ in the remaining domain,
until it exits $D_1$ 
 or its $\phi_0$-capacity exceeds $\eta_1$.
\end{enumerate}
\end{Lem}

\begin{proof}
Informally, the argument is the following: divide the two SLEs in $n$
segments; one has to prove that one can either grow the $n$ segments
of the first SLE, then the $n$ segments of the second SLE, or the
other way round and get the same law. The permutation of two segments
(of the two SLEs) induces an error term of $O(n^{-3})$, from the
infinitesimal commutation relations. One needs $n^2$ such
permutations; letting $n$ go to infinity, one gets the result. The
uniformity in the error terms is provided by the restriction to paths
in the disjoint compact neighbourhoods $D_1$, $D_2$.

For simplicity, we will consider only the case where $\phi_0=\Id$ (and
the $\phi_0$-capacity is the ordinary half-plane capacity). For the general case,
one has to replace fixed times by corresponding stopping times; the
proof goes otherwise unchanged.

For positive times $S$, $T$, let $\E_1$ designate the expectation for
pairs of random curves obtained by growing first the $\SLE$ in $D_1$
up to time $S$ (half-plane capacity $2S$), and then the $\SLE$ in $D_2$
up to time $T$. The symbol $\E_2$ refers to expectation for the
reversed construction. The driving process for each of these Loewner
chains (seen in the original half-plane) is denoted by $X$, $Y$. Let
$\tau_1$ be the time at which the first SLE exits $D_1$, and $\tau_2$
the corresponding time for the second SLE. We will prove that
$$\E_1\left(.\ind_{\tau_1>S,\tau_2>T}\right)=\E_2\left(.\ind_{\tau_1>S,\tau_2>T}
\right)$$
as measures on $C_0([0,S],\R)\times C_0([0,T],\R)$. To recover the
statement of the lemma, one then considers the measures:
$$\E_i\left(.(\ind_{\tau_1>S,\tau_2>T}-\ind_{\tau_1>S+ds,\tau_2>T})\right)$$
for $i=1,2$. So we can work with fixed times $S$ and $T$. Note that
$2\tau_i$ is bounded by the half-plane capacity of $D_i$, $i=1,2$.

Let $S_0=0<S_1<\cdots< S_m=S$ and $T_0=0<T_1<\cdots< T_m=T$ be
fixed sequences of times ($m\geq 1$). Also, let $(\tilde\varphi_i)_{0\leq
i\leq m}$, $(\tilde\psi_i)_{0\leq i\leq m}$ be test functions (i.e. in
$C^\infty_c(\R)$). By a monotone class argument, we need only to see that:
$$\E_1\left(\prod_k\tilde\varphi_k(X_{S_k})\tilde\psi_k(Y_{T_k})\ind_{\tau_1>S,\tau_2>T}\right)=
\E_2\left(\prod_k\tilde\varphi_k(X_{S_k})\tilde\psi_k(Y_{T_k})\ind_{\tau_1>S,\tau_2>T}\right)$$

For $n\geq 1$, consider increasing
sequences $(s_i)_{0\leq i\leq mn}$, $(t_i)_{0\leq i\leq mn}$, where
$s_{nj}=S_j$, $t_{nj}=T_j$, and the increments $(s_{i+1}-s_i)$,
$(t_{i+1}-t_i)$ go uniformly to 0. Define $\varphi_{ni}=\tilde\varphi_i$,
$\psi_{ni}=\tilde\psi_i$, and $\varphi_i=\psi_i=1$ if $n$ does not
divide $i$.

Note that the commutation relation holds for functions of the positions
of all marked points in the Loewner flow. For convenience, we will
approximate the event $\{\tau_1>S,\tau_2>T\}$ by a function of an
extended flow. More precisely, let $\delta>0$ be a (small) positive
number and $N$ a (large) integer. Mark $N$ points $z_1,\dots z_N$ 
on the Jordan boundaries of $D_1$ and $D_2$ (one can also mark their
conjugates $\overline z_1,\dots,\overline z_N$, extending the flow by
Schwarz reflection). For instance, one can
choose them so that the Hausdorff distance between $\partial
D_1\cup\partial D_2$ and $\{z_1,\dots z_N\}$ is minimal ($N$ being
fixed).

Let ${\mc K}=\{(K,u)\}$ where $K$ is a compact hull included in $D_1$
and $u$ is a point in $\partial D_1\cap K$. Then  ${\mc K}$ is a
compact set (using the Hausdorff metrics on compact subsets of
$D_1$), so one can choose $\delta>0$ so that:
$$\delta>\max_{(K,u)\in{\mc K}}\left(\min_i |\phi_K(z_i)-\phi_K(u)|\right)$$
and the corresponding inequality holds for hulls in $D_2$. Also, it is
easy to see that one can choose $\delta$ so that it goes to zero as $N$
goes to infinity, by a compacity argument.
Let $f_\delta$ be a smooth function of the
variables $(x,y,z_1,\dots z_N)$, taking values in $[0,1]$, such that it vanishes if $|x-z_i|<\delta$
or $|y-z_i|<\delta$ for some $i\in\{1,\dots,N\}$ and equals 1 if
$|x-z_i|>2\delta$, $|y-z_i|>2\delta$ for all $i$. One can assume that one of
the $z_i$'s is real and between $x$ and $y$, and similarly, the 
other marked points (that influence the drift) are separated from $x$,
$y$ by one of the ``spectator'' $z_i$'s. The choice of
$\delta$ ensures that $f_\delta$ vanishes as soon as an SLE crosses $\partial
D_i$, $i=1,2$.

Let $W$ denote the full configuration (images of growth points and marked
points in the Loewner flow). Then we just have to prove that:
$$\lim_{n\rightarrow\infty}(\E_1-\E_2)\left(\prod_k\varphi_k(X_{s_k})f_\delta(W_{(s_k,0)})\psi_k(Y_{t_k})f_\delta(W_{(0,t_k)})\right)=0$$
and then let $N\nearrow\infty$, $\delta\searrow 0$ to get the result for
stopped SLEs. So it what follows we may replace
$\varphi_k(X_{s_k})f_\delta(W_{(s_k,0)})$ with
$\varphi_k(W_{(s_k,0)})$ (a function of the configuration), and
similarly $\psi_k(X_{s_k})f_\delta(W_{(s_k,0)})$ with
$\psi_k(W_{(s_k,0)})$

Consider also two random curves $\gamma$, $\hat\gamma$ started from
$x$ (resp. $y$) in $\H$, parameterized by half-plane capacity.
%in $\H$, and two random curves $\hat\gamma$, $\hat\gamma'$
%started from $y$, such that
%$\gamma\cap\hat\gamma=\gamma'\cap\hat\gamma'=\varnothing$; $h,h',\hat
%h,\hat h'$ are the associated conformal equivalences,
%$x_i=h_{s_i}(\gamma_{s_i}),\dots, y_i=\hat h_{t_i}(\hat\gamma_{t_i})$ (we
%shall use the same notations for $\gamma'$, $\hat\gamma'$). 
Let $\phi_{s,t}=\phi_{\gamma_{[0,s]}\cup\hat\gamma_{[0,t]}}$, and $W_{s,t}$
is the configuration $\phi_{s,t}(\gamma_s,\hat\gamma_t,z_1,\dots
z_N,\dots)$. We will also abbreviate $x_k=W_{(s_k,0)}$, $y_k=W_{(0,t_k)}$.

Consider two permutations $\sigma$ and $\sigma'$ of $\{s_1,\dots
s_{mn},t_1,\dots t_{mn}\}$, increasing for the partial order generated
by $s_k< s_{k+1}$, $t_k<t_{k+1}$, and such that $\sigma$ and
$\sigma'$ differ by a transposition of two consecutive elements. For
instance $\sigma=(s_1,\dots s_{mn},t_1,\dots t_{mn})$ and
$\sigma'=(s_1,\dots s_{mn-1},t_1,s_{mn},t_2,\dots t_{mn})$. 
Suppose that $\gamma,\hat\gamma$ are obtained from the permutation $\sigma$ in the following fashion:
if $\sigma=(\sigma_1, s_{i+1}, \sigma_2)$, $t_j$ is the maximal $t_.$
element in $\sigma_1$, and
$\phi=\phi_{s_i,t_j}$, then
$\phi(\gamma_{[s_i,s_{i+1}]})$ is an $\SLE_\kappa(b)$ started from
$W_{s_i,t_j}$ and independent of
$\phi$ conditionally on its starting state
(stopped so that $\phi_{s_i,0}(\gamma_{[s_i,s_{i+1}]})$ has capacity $2(s_{i+1}-s_i)$). Likewise, if $\sigma=(\sigma_1, t_{j+1}, \sigma_2)$, $s_i$
is the maximal $s_.$ element in $\sigma_1$, and
$\phi=\phi_{s_i,t_j}$, then
$\phi(\gamma_{[t_j,t_{j+1}]})$ is an $\SLE_{\tilde\kappa}(\tilde b)$ started from
$W_{s_i,t_j}$ and independent of
$\phi$ conditionally on its starting state. The symbol $\E$ is
expectation for this construction (relative to $\sigma$), and $\E'$ is
the corresponding expectation obtained from $\sigma'$.

Let $\sigma=(\sigma_1,s_i,t_j,\sigma_2)$ and
$\sigma'=(\sigma_1,t_j,s_i,\sigma_2)$. Then:
\begin{align*}
\E'\left(\prod_k\varphi_k(x_k)\psi_k(y_k)\right)&=\E'\left(\left(\prod_{k<i,l<j}\varphi_k(x_k)\psi_l(y_l)\right)\varphi_i(x_i)\psi_j(y_j)\left(\prod_{k>i,l>j}\varphi_k(x_k)\psi_l(y_l)\right)\right)
\end{align*}
(Here the $f_\delta$ are implicitly included in the $\varphi_k$, $\psi_k$).
The expectation of the last part of the product conditionally on
$\gamma_{[0,s_i]},\hat\gamma_{[0,t_j]}$ is a function of
$W_{s_i,t_j}$,
Denote
by $F(u,v)$ this function, which is the same under $\E$ and $\E'$; by
induction and standard regularity results (the drift terms stay
bounded as long as the functional does not vanish), it is
easily seen that $F$ is a smooth function ;
the existence of regular conditional probability is clear for
the same reasons.

Now, consider:
$$\E\left(\left.
\varphi_i(x_i)\psi_j(y_j)F(\phi_{s_i,t_j}(\gamma_{s_i},\hat\gamma_{t_j}))\right|W_{s_{i-1},t_{j-1}}\right)-\E'\left(\left.
\varphi_i(x_i)\psi_j(y_j)F(\phi_{s_i,t_j}(\gamma'_{s_i},\hat\gamma'_{t_j}))\right|W_{s_{i-1},t_{j-1}}\right)$$
Assume that $i,j$ are not multiples of $n$ (and traces are away from
the boundaries of $D_1$, $D_2$ at time $i,j$). Then $\varphi_i=\psi_j=1$, and this difference is
$O(n^{-3})$, from the infinitesimal commutation relation. 

If $i$ or $j$ is a multiple of $n$, note that, if $x'=\phi_{s_i,t_{j-1}}(\gamma_{s_i})$, $y'=\phi_{s_{i-1},t_{j}}(\gamma_{t_j})$, 
 $x''=\phi_{s_i,t_{j}}(\gamma_{s_i})$,
$y''=\phi_{s_{i},t_{j}}(\gamma_{t_j})$, then:
$$x'=x''-\frac {2(t_j-t_{j-1})}{x''-y''}+O(n^{-2})$$
using the backward Loewner flow; one gets a similar expression for
$y'$, and these hold under $\E$ and $\E'$. So in this case the
difference is $O(n^{-2})$. 

To get from $\sigma=(s_1,\dots,s_{mn},s_1,\dots,s_{mn})$ to
$\sigma'=(t_1,\dots,t_{mn},t_1,\dots,t_{mn})$, one needs $(mn)^2$
transpositions ($(mn)$ transpositions to bring $t_1$ in first
position, then $(mn)$ transpositions to bring $t_2$ in second
position, ...). For such a transposition $(s_i,t_j)$, $i$ or $j$ is a
multiple of $n$ in $m^2(2n-1)$ case. This transposition is valid as
long as the paths stay in $D_1,D_2$ (more precisely, as long as the
$f_\delta$ terms are 1). Conversely, if a path is close to the
boundary of $D_1,D_2$, the functional is zero with probability close
to one. Hence: 
%$$\E\left(\prod_{k=1}^m\tilde\varphi_k(\phi_{S_k,0}(\gamma_{S_k}))\tilde\psi_k(\phi_{0,T_k}(\hat\gamma_{T_k}))\right)-\E'\left(\prod_{k=1}^m\tilde\varphi_k(\phi_{S_k,0}(\gamma'_{S_k}))\tilde\psi_k(\phi_{0,T_k}(\hat\gamma'_{T_k}))\right)=m^2(n-1)^2O(n^{-3})+m^2(2n-1)O(n^{-2})$$
%$$(\E_1-\E_2)\left(\prod_k\varphi_k(X_{s_k})f_\delta(W_{(s_k,0)})\psi_k(Y_{t_k})f_\delta(W_{(0,t_k)})\right)=m^2(n-1)^2O(n^{-3})+m^2(2n-1)O(n^{-2})$$
$$(\E_1-\E_2)\left(\prod_k\varphi_k(X_{s_k})\psi_k(Y_{t_k})\ind_{E(n,N,\delta)}\right)=m^2(n-1)^2O(n^{-3})+m^2(2n-1)O(n^{-2})=O(C(N,\delta)/n)$$
where $X_s=\phi_{s,0}(\gamma_s)$, $Y_t=\phi_{0,t}(\hat\gamma_t)$, and
$E(n,N,\delta)$ is the event than none of the $f_\delta$'s vanishes at
a sampled time. The error term is uniform in $n$ but depends on
$N,\delta$.

As $n$ goes to infinity
($N$, $\delta$ being fixed), the
probability that the first SLE crosses $\partial D_1$ without $f_\delta$
vanishing at one of the sampled times $s_i$ goes to zero (since in
this case $f_\delta$ vanishes on an open set of times). So we can
assume that the SLEs stay in $D_1$, $D_2$, hence we have uniformity in
the $O(n^{-3})$ estimate of the commutation condition. The last case
to study is when the trace gets close to the boundary,
say $|x_{i-1}-z_j|<2\delta$ for some $i,j$, without actually crossing
it. The probability of this event goes to zero as $N$ goes to infinity
and $\delta$ goes to zero.

So the above estimate is valid up to an event of negligible
probability, viz. either an SLE crosses $\partial D_1$ or $\partial D_2$
without the functional vanishing or one of the $f_\delta$ is less than
one and yet the functional does not vanish.
Taking the limit as $n$ goes to infinity and 
$N\nearrow\infty$, $\delta\searrow 0$ (so that
$C(N,\delta)/n\rightarrow 0$), one gets the stated
identity, that is :
$$\E_1\left(.\ind_{\tau_1>S,\tau_2>T}\right)=\E_2\left(.\ind_{\tau_1>S,\tau_2>T}
\right).$$
This concludes the proof of the lemma.

%Note that we have implicitly assumed uniformity in these
%estimates. This can be justified as follows: in this situation
%(commutation of two $\SLE_\kappa(2)$, $\kappa\leq 4$), the curves
%$\gamma,\hat\gamma$ (resp. $\gamma',\hat\gamma'$) are disjoint. As a
%consequence,
%$(\phi_{s,t}(\gamma_{[0,s]})-\phi_{s,t}(\hat\gamma_{[0,t]}))_{s\leq
%S,t\leq T}$ is a.s. uniformly bounded away from 0, which is enough to
%ensure uniformity in the commutation relation.   

\end{proof}

\begin{proof}[Proof of the Proposition.]
Let $D_1$, $D_2$ be as in the lemma. We grow an $\SLE_\kappa(b)$ in $D_1$ until it
reaches $\partial D_1$, and then in the remaining domain an
$\SLE_{\tilde\kappa}(\tilde b)$ in $D_2$ until it
reaches $\partial D_2$. This defines a Loewner chain
$(K_{s,t})_{(s,t)\in{\mc T}}$. We will prove that this chain is an
$\SLE(\kappa,b,\tilde\kappa,\tilde b)$, i.e. it has the appropriate
Markov property.

Let $0=S_0<S_1<\cdots<S_k=\infty$ and $0=T_0<T_1<\cdots<T_k=\infty$ be sequences
of fixed times. Let $\sigma$ be a permutation of the symbols
$(S_1,\dots, S_k, T_1,\dots, T_k)$,  which is
increasing for the partial order generated by $S_i<S_{i+1}$,
$T_i<T_{i+1}$. Let $(K^\sigma_{s,t})$ be the (random) Loewner chain obtained by
growing SLEs alternatively in $D_1$ and in $D_2$ according to
$\sigma$, stopping the SLEs when they reach $\partial D_1,\partial
D_2$. For instance, if $\sigma=(S_1,T_1,T_2,S_2,\dots)$, one grows the first
$\SLE$ to half-plane capacity $2S_1$ (and stop it if it reaches $\partial
D_1$), then the second $\SLE$ to half-plane capacity $2S_2$, measured
in the original half-plane (and stop it if it reaches $\partial
D_1$), and then again the first
$\SLE$ to half-plane capacity $2S_2$ (and stop it \dots), \dots. 
To alleviate notations, we will use the convention that for a Loewner
chain $(\tilde K_{s,t})$, $\tilde K_{s,t}=\tilde K_{s\wedge S,t\wedge T}$, where $S$ is the
time at which $(\tilde K_{s,0})_s$ exits $D_1$ (resp.  $T$ is the
time at which $(\tilde K_{0,t})_t$ exits $D_2$).

If $\sigma$ and $\sigma'$ differ by a single tranposition,
i.e. $\sigma=(\sigma_1,S_{k_1},T_{k_2},\sigma_2)$,
$\sigma'=(\sigma_1,T_{k_2},S_{k_1},\sigma_2)$, then we can apply the
lemma with $\phi_0=\phi_{K^{\sigma_1}_{S_{k_1},T_{k_2}}}$, $D'_1=\phi_0(D_1)
$, $D'_2=\phi_0(D_2)$. This proves that we can couple
$(K^\sigma_{s,t})$ and $(K^{\sigma'}_{s,t})$. By induction, we can
couple (simultaneously) the chains $(K^\sigma_{s,t})$ for all
admissible permutations $\sigma$. By construction, for
$\sigma=(S_1,\dots S_k, T_1,\dots, T_k)$, $K^\sigma$ is distributed as $K$.

Now, for any $k_1,k_2\in\{0,\dots k\}$, one can consider the
permutation:
$$\sigma=(S_0,\dots,S_{k_1},T_0,\dots,T_{k_2},S_{k_1+1},\dots,S_k,T_{k_2+1},\dots
T_k).$$ 
The previous coupling proves the Markov property for the fixed
time $(S_{k_1},T_{k_2})$ (i.e the chain
$$(g_{S_{k_1},T_{k_2}}(K_{S_1+s,T_2}\setminus K_{S_1,T_2}))$$
 is a stopped $\SLE_\kappa(b)$, and the same thing holds for the other $\SLE$).

This still holds for stopping times supported on
$\{(S_{k_1},T_{k_2}),k_1,k_2=0,\dots k\}$. Since the subdivisions
$S_0<\cdots< S_k$ and $T_0<\cdots< T_k$ were arbitrary, this also holds
for stopping times with finite support, and by a limiting argument for
all stopping times (as for the classical Markov property). 
 
\end{proof}

\section{Classification of commuting SLEs}

We can now conclude the general study of pairs of commuting chordal SLEs in a
simply connected domain. In the upper half-plane $\H$, with $(2n+2)$
marked points $(x,y,z_1,\dots,z_n)$ on the real line (and one marked
point $z_{n+1}$ at infinity), consider two
parameters $\kappa,\tilde\kappa$, and two smooth functions of the
configuration $(x,y,z_1,\dots,z_n)$,
translation invariant and homogeneous of degree $-1$. Let ${\mc L}$ be
the infinitesimal generator of the $\SLE_\kappa(b)$ growing at $x$
(driven by $(X_s)$, $(g_s)$ are the corresponding conformal
equivalences), and
${\mc M}$ 
the infinitesimal generator of the $\SLE_{\tilde\kappa}(\tilde b)$
growing at $y$ (driven by $(Y_t)$, $(\tilde g_t)$ are the corresponding conformal
equivalences). By a cocycle $C(\phi,z)$ we mean a function of the form:
$$C(\phi,z)=\prod_{1\leq i<j\leq n+1}\left(\phi'(z_i)\phi'(z_j)\left(\frac{z_j-z_i}{\phi(z_j)-\phi(z_i)}\right)^2\right)^{\nu_{ij}}\frac{f(\phi(z_1),\dots,\phi(z_n))}{f(z_1,\dots,z_n)}$$
where $f$ is a non-vanishing function (translation invariant and
homogeneous of degree 0). So $f$ is a conformally invariant function
of the marked points $(z-1,\dots,z_n,z_{n+1}=\infty)$ and can be seen
as a function on the residual (i.e. not invovlving the positions of
$x,y$) moduli space.
If $n\geq 3$, this decomposition is not
unique, as discussed in Section 4.

\begin{Thm}
The following assertions are equivalent.
\begin{enumerate}
\item There exists an $\SLE(\kappa,b,\tilde\kappa,\tilde b)$.
\item The infinitesimal generators satisfy the relation:
\begin{equation*}
\left[{\mc L},{\mc M}\right]=\frac{4}{(y-x)^2}\left({\mc M}-{\mc L}\right)
\end{equation*}
\item $\tilde\kappa=\kappa$ or $\tilde\kappa=16/\kappa$,
$b=\kappa\partial_x\psi/\psi$, $\tilde
b=\tilde\kappa\partial_y\psi/\psi$, where $\psi$ is
a non-vanishing solution of the system:
$$\left\{\begin{array}{l}
\displaystyle\frac{\kappa}
2\partial_{xx}\psi+\sum_i\frac{2\partial_i\psi}{z_i-x}+\frac{2\partial_y\psi}{y-x}+\left(\left(1-\frac
6{\tilde\kappa}\right)\frac 1{(y-x)^2}+h(x,z)\right)\psi=0\\
\displaystyle\frac{\tilde\kappa}2\partial_{yy}\psi+\sum_i\frac{2\partial_i\psi}{z_i-y}+\frac{2\partial_x\psi}{x-y}+\left(\left(1-\frac
6\kappa\right)\frac 1{(x-y)^2}+h(y,z)\right)\psi=0
\end{array}\right.$$
where 
$$h(x,z)=\sum_i\frac{\mu_i}{(x-z_i)^2}+\sum_{i<j}\frac{\nu_{ij}}{(x-z_i)(x-z_j)}+\ell_x\log(f(z)),$$
the $\mu_i$, $\nu_{ij}$ are constant parameters, and $f$ is a function
on the residual moduli space.

\item $\tilde\kappa=\kappa$ or $\tilde\kappa=16/\kappa$,
and there is 
a non-vanishing function $\psi$ and a cocycle $C$ such that if:
\begin{align*}
Z_s=&\psi(X_s,g_s(y),\dots,g_s(z_i),\dots)g'_s(y)^{\alpha_{\tilde\kappa}}C(g_s,z)\\
\tilde Z_t=&\psi(\tilde g_t(x),Y_t,\dots,\tilde g_t(z_i),\dots)\tilde
g'_t(x)^{\alpha_\kappa}
C(\tilde g_t,z)
\end{align*}
then $(Z_s)$ is a local martingale for chordal
$\SLE_\kappa(x\rightarrow\infty)$, $(\tilde Z_t)$ is a local martingale for chordal
$\SLE_{\tilde\kappa}(y\rightarrow\infty)$, and $\SLE_\kappa(b)$ is the
Girsanov transform of chordal $\SLE_\kappa$ by $(Z_s)$,
$\SLE_{\tilde\kappa}(\tilde b)$ is the
Girsanov transform of chordal $\SLE_{\tilde\kappa}$ by $(\tilde Z_t)$. 
\end{enumerate}
\end{Thm}

Note that there is no loss of generality in considering two (rather
than $m\geq 2$) commuting SLEs. Indeed, the only conditions will be
the pairwise conditions. It is also easy to see that the proofs for
local commutation can be adapted for $m$ $\SLE$s (though notations
become quite heavy). Let us explicit, say, condition (iii) in this
situation. On the real line, $(m+n)$ points $(y_1,\dots, y_m,z_1,\dots,
z_n)$ are marked, and we want to grow $m$ SLEs ($\SLE_{\kappa_i}(b_i)$,
$i=1,\dots, m$) at $y_1,\dots,y_m$. Then
$\{\kappa_1,\dots,\kappa_m\}\subset\{\kappa_1,16/\kappa_1\}$,
$b_i=\kappa_i\partial_{y_i}\psi/\psi$, where $\psi$ is annihilated by
the operators:
$$\frac{\kappa}
2\partial_{y_iy_i}+\sum_{j\neq i}\frac{2\partial_{y_j}}{y_j-y_i}+
\sum_{j}\frac{2\partial_{z_j}}{z_j-y_i}
-2\left(\sum_{j\neq i}\frac{\alpha_{\kappa_j}}{(y_j-y_i)^2}+\sum_j\frac{\mu'_j}{(z_j-y_i)^2}+\sum_{j<j'}\frac{\nu'_{jj'}}{(z_j-y_i)(z_{j'}-y_i)}+\ell_{y_i}\log(f(z))
\right)$$
for some parameters $\mu'_j$, $\nu'_{jj'}$ and some function $f$ on
the residual moduli space.

\section{Restriction formulae for non-intersecting SLEs}

In this section we specialize to a simple parametric case, where $n$
SLEs started from distinct points on the real line are aiming at
infinity; there are only $n$ marked points on the real line (and one
at infinity). Each of the $n$ SLEs is an
$\SLE_\kappa(\underline\rho)$, where $\underline\rho=(2,\dots 2)$. 
In this situation, we can not only define locally a $n$-parameter
Loewner chain, but also define it globally if $\kappa\leq 4$. Indeed,
the only thing preventing from a global definition is the possibility
of collisions of marked points. But such collisions a.s. don't happen
for these $\SLE_\kappa(2,\dots,2)$, so we can actually define a chain
with full time set $\R_+^n$. 

If the $n$ starting points collapse to zero, we get $n$
``non-intersecting'' SLEs starting at $0$ and ending at
$\infty$. Restriction formulae are derived for these Loewner chains
(indexed by $\R_+^n$). This gives a simple realization of the
exponents $h_{1;n+1}(\kappa)$ (see also \cite{W2}).

The radial case ($n$ ``non-intersecting'' SLEs started from the
boundary and aiming at a single bulk point) is also studied, and
restriction formulae then give the exponent $2h_{0;n/2}(\kappa)$.

\subsection{The chordal case}

Let $y_1<\cdots<y_n$ be $n$ real points. Consider the infinitesimal generators:
$${\mc L}_i=\frac\kappa 2\partial_{ii}+\left(\sum_{j\neq i}\frac
{2}{y_i-y_j}\right)\partial_i+\sum_{j\neq i}\frac
2{y_j-y_i}\partial_j$$
Then, from the previous computations (parametric case), we see that the following
commutation relations are satisfied:
$$\left[{\mc L}_i,{\mc L}_j\right]=\frac{4}{(y_j-y_i)^2}\left({\mc
L}_j-{\mc L}_i\right).$$

As mentioned earlier, this ensures (if $\kappa\leq 4$) the
existence of a process $(K_{s_1,\dots,s_n})$ such that:\\
$(\phi_{K_{s_1,\dots,s_n}}(K_{s_1,\dots s_i+s, \dots s_n}))_{s\geq 0}$
is an $\SLE_\kappa(2,\dots 2)$ started from
$(\phi(\gamma^{i}_{s_i}))=(\phi(\gamma^{1}_{s_1}),\dots,\phi(\gamma^{n}_{s_n}))$,
independent from $(K_{t_1,\dots t_n})_{t_j\leq s_j}$, where
$\phi=\phi_{K_{s_1,\dots,s_n}}$, $K_{0,\dots, s_i,\dots
0}=\gamma^i_{[0,s_i]}$. Assume that this process is started from
$(y_1,\dots,y_n)$, i.e. $\gamma^i_0=y_i$, and define
$K_\infty=\bigcup_i\gamma^i_{[0,\infty]}\subset\overline\H$.

Consider now a hull $A\subset\overline\H$, that does not intersect
$\{y_1,\dots,y_n\}$.  
Let $\lambda_\kappa=(6-\kappa)(8-3\kappa)/2\kappa$. If $\kappa\leq
8/3$, and $L$ is an independent random loop soup with intensity
$\lambda_\kappa$ in $\H$, define
$K_\infty^L$ to be (the filling of) the union of $K_\infty$ and the
loops in $L$ that intersect it. 
Then, if $\kappa\leq 8/3$, the following restriction formula holds:

\begin{Lem} The probability that $K_\infty^L$ does not intersect
$A$ is given by: 
$$\P(K_\infty^L \cap
A=\varnothing)=\prod_i\phi_A'(y_i)^{(6-\kappa)/2\kappa}\prod_{i<j}\left(\frac{\phi_A(y_j)-\phi_A(y_i)}{y_j-y_i}\right)^{2/\kappa}$$
\end{Lem}
\begin{proof}
Define ${\underline s}=(s_1,\dots s_n)$, $h_{\underline
s}=\phi_{\phi_{K_{\underline s}}(A)}$, and $Y^{(i)}_{\underline
s}=\phi_{\underline s}(\gamma^i_{s_i})$. Then, from the definition of
$K_{\underline s}$ and Lemma 4 in \cite{D4}, one sees that: 
$$s_k\longmapsto \prod_i h'_{\underline s}(Y^{(i)}_{\underline
s})^{(6-\kappa)/2\kappa}\prod_{i<j}\left(\frac{h_{\underline s}(y_j)-
h_{\underline
s}(y_i)}{y_j-y_i}\right)^{2/\kappa}\exp\left(\lambda_\kappa\int_{0}^{s_k}\frac{Sh_{s_1,\dots
\sigma_k, \dots s_n}(Y^{(k)}_{s_1,\dots \sigma_k, \dots
s_n})}6d\sigma_k\right)$$ 
is a bounded martingale for all $k$, $s_1,\dots, s_n$. From the
properties of the Brownian loop soup, it appears that the following
semimartingale (proportional to the first one) is also a bounded
martingale: 
$$s_k\longmapsto \prod_i h'_{\underline s}(Y^{(i)}_{\underline
s})^{(6-\kappa)/2\kappa}\prod_{i<j}\left(\frac{h_{\underline s}(y_j)-
h_{\underline s}(y_i)}{y_j-y_i}\right)^{2/\kappa}\P(K_{\underline
s}^L\cap A=\varnothing)$$ 
Hence, for all ${\underline s}=(s_1,\dots s_n)$, one gets (using $n$
different martingales): 
$$\prod_i\phi_A'(y_i)^{(6-\kappa)/2\kappa}\prod_{i<j}\left(\frac{\phi_A(y_j)-\phi_A(y_i)}{y_j-y_i}\right)^{2/\kappa}=\E\left(\prod_i
h'_{\underline s}(Y^{(i)}_{\underline
s})^{(6-\kappa)/2\kappa}\prod_{i<j}\left(\frac{h_{\underline s}(y_j)-
h_{\underline s}(y_i)}{y_j-y_i}\right)^{2/\kappa}\P(K_{\underline
s}^L\cap A=\varnothing)\right)$$ 
Now, as $\inf\underline s$ goes to infinity, the product in the
right-hand side converges to $\P(K_\infty^L \cap A=\varnothing)$,
which concludes the proof. 
\end{proof}

Define the conformal weight $h_{p;q}=h_{p;q}(\kappa)$ by:
$$h_{p;q}=\frac{(p\kappa-4q)^2-(\kappa-4)^2}{16\kappa}$$
Then, if $y_1,\dots y_n$ collapse to zero, the above formula reduces to:
$$\P(K_\infty^L \cap A=\varnothing)=\phi_A'(0)^{h_{1;n+1}}$$
The role of the conformal weights $h_{1;n+1}$ in the context of restriction
measures and $\SLE_\kappa(\rho)$ is discussed in \cite{W2}.

\begin{Cor}[Restriction property]
Let $\kappa\leq 8/3$, $L$ a loop soup with intensity $\lambda_\kappa$.
Conditionally on $\{K_\infty^L\cap A=\varnothing\}$
and up to a time change, $((\phi_A(K_{\underline s}))_{\underline
s},\phi_A(\{\delta\in L:\delta\cap A=\varnothing\}))$
is distributed as $((K_{\underline s})_{\underline s},L)$, where
$(K_{\underline s})$ is started from $\phi_A(y_1),\dots
\phi_A(y_n)$. In particular, the collection of measures on hulls
$K_\infty^L$ indexed by the starting configuration
$(\H,z_1,\dots,z_n,z_{n+1}=\infty)$ has the restriction property.
\end{Cor}
\begin{proof}
From the previous lemma, the result is a straightforward application
of the Girsanov theorem and the restriction property of the loop soup.
The assertion on the restriction property can be derived directly by
applying the previous formula to concatenation of hulls $A.B$. 
\end{proof}

\subsection{The radial case}

Recall the definition of radial $\SLE_\kappa(\underline\rho)$:
assume the existence of processes $(\xi_t)_{t\geq 0}$ and
$(\chi^{(i)}_t)_{t\geq 0}$, $i\in\{1\dots k\}$, satisfying the SDEs:

\begin{equation*}
\left\{\begin{array}{l}d\xi_t=(i\xi_t\sqrt\kappa dB_t-\frac\kappa
2\xi_tdt)+\sum_{i=1}^k\frac{\rho_i}2\left(-\xi_t\frac{\xi_t+\chi^{(i)}_t}{\xi_t-\chi^{(i)}_t}\right)dt\\
d\chi^{(i)}_t=-\chi_t\frac{\chi_t+\xi_t}{\chi_t-\xi_t}dt\end{array}\right.
\end{equation*}

Then the ODEs $dg_t(z)=-g_t(z)(g_t(z)+\xi_t)/(g_t(z)-\xi_t)dt$ define
radial $\SLE_\kappa(\underline\rho)$ in the unit disk $\U$.  

First, we briefly discuss commutation conditions in the radial
case. Suppose that $\chi_1,\dots,\chi_n$ are $n$ points on
the unit circle. One considers two SLEs growing at $\chi_1$ and
$\chi_n$ resp., assuming that the drift terms are functions of the
$\chi_i$. We think of functions annihilated by infinitesimal
generators as expected values of some events; it is quite natural to
express these real-valued functions in angular coordinates:
$\chi_j=\exp(i\theta_j)$. Reasoning as in the chordal case,
if ${\mc L}$ and ${\mc M}$ are the infinitesimal generators of the two
SLEs, the commutation condition reads:
$$[{\mc L},{\mc M}]=\frac 1{\sin((\theta_n-\theta_1)/2)^2}({\mc L}-{\mc M}).$$
The generator for a radial $\SLE_\kappa(\rho_1,\dots,\rho_{n-1})$
started from $(\chi_1,\dots,\chi_n)$ is:
$${\mc L}=\frac\kappa 2\partial_{11}+\sum_{i>1}\cot\left(\frac{\theta_i-\theta_1}2\right)\left(\partial_i-\frac{\rho_i}2\partial_1\right)$$
By analogy with the chordal case, one can find solutions for this
commutation relation:
\begin{enumerate}
\item 
Two $\SLE_\kappa(\rho)$ started from $(\chi_1,\chi_2)$ and
$(\chi_2,\chi_1)$ resp., $\rho\in\{2,\kappa-6\}$.
\item $n$ $\SLE_\kappa(2,\dots,2)$ started from
$(\chi_i,\chi_1,\dots,\widehat{\chi_i},\dots\chi_n)$.
\end{enumerate}

Let us comment briefly on the case (i). For $\rho=\kappa-6$, this is
only chordal reversibility in a radial normalization (as may be seen
by slightly modifying the argument for chordal-radial equivalence when
$\kappa=6$). In the case $\rho=2$, this gives a model of ``pinned
chordal SLE'', i.e chordal SLE ``conditionally'' on the trace visiting
a given bulk point, for $\kappa<8$. More precisely, start from a
chordal SLE in radial normalization (hence, up to a time change,
radial $\SLE_\kappa(\kappa-6)$). Then the first moment estimate in
\cite{B1} relies on the computation of the leading eigenvector for the
associated infinitesimal generator. This yields a local martingale:
$$s\mapsto\sin\left(\frac{\theta^{(1)}_s-\theta^{(2)}_s}2\right)^{8/\kappa-1}e^{(1-\kappa/8)s}$$
corresponding of the probability that the SLE trace gets infinitely
close to the bulk point 0. Using this as a Girsanov density, one gets
a radial $\SLE_\kappa(2)$. Note that for $\kappa=8$, this density is
1, and $\kappa-6=2$.\\
There are other examples with two boundary points. Consider a chordal
$\SLE_\kappa$ from $\chi_1$ to $\chi_2$ and condition it to leave $0$
on its left (resp. right); this can be made explicit (see
\cite{S1}). Once again, the drift terms are (generically in $\kappa$) transcendental.

In the case (ii), just as in the chordal case (at least if $\kappa\leq
4$), based on the infinitesimal commutation relations, one can define
a $n$-braids radial SLE. The question of such a definition, from a CFT
point of view, appears in \cite{Ca3,Ca3c}. 
Note also that summing the $n$ generators here gives the generator of 
Dyson's Brownian motion.

As above, we study the case (ii) from the restriction point of view.
First, we have to derive restriction formulae for radial $\SLE_\kappa(\underline\rho)$.
Let $A$ be a hull of $\overline \U$ (i.e. A is a compact subset of $\overline \U$,
$A\cap\overline\U=\overline{(A\cap\U)}$, $\U\setminus A$ is equivalent
to $\U$ and $0\notin A$). For any such hull, denote by $\phi_A$ the
unique conformal equivalence $\U\setminus A\rightarrow \U$ such that
$\phi_A(0)=0$ and $\phi_A'(0)>0$. Suppose that $\xi_0,\chi^{(i)}_0$
are not in $A$. Then $h_t=\phi_{g_t(A)}$ is defined at least for small
times. Recall that
$\lambda_\kappa=(8-3\kappa)(6-\kappa)/2\kappa$. Then the following
result (analogous to Lemma 4 in \cite{D4} and generalizing a result
stated in \cite{LSW3}) holds: 

\begin{Lem}
In this situation, define $\tilde\xi_t=h_t(\xi_t)$, $\tilde\chi^{(i)}_t=h_t(\chi^{(i)}_t)$:
\begin{align*}
M^\varnothing_t&=\left(\frac{h'_t(\xi_t)\xi_t}{\tilde\xi_t}\right)^{\frac{(6-\kappa)}{2\kappa}} h'_t(0)^{\frac{(6-\kappa)(\kappa-2)}{8\kappa}}\exp\left(-\lambda_\kappa\int_0^t\frac{\xi_s^2Sh_s(\xi_s)}6ds\right)\\
M^{(i)}_t&=\left(\frac{h'_t(\chi^{(i)}_t)\chi^{(i)}_t}{\tilde\chi^{(i)}_t}\right)^{\frac{\rho_i(\rho_i+4-\kappa)}{4\kappa}}\left(\frac{(\tilde\chi^{(i)}_t-\tilde\xi_t)^2}{\tilde\chi^{(i)}_t\tilde\xi_t}\frac{\chi^{(i)}_t\xi_t}{(\chi^{(i)}_t-\xi
_t)^2} \right)^{\frac{\rho_i}{2\kappa}}h'_t(0)^{\frac{\rho_i(\rho_i+4)}{8\kappa}}\\
M^{(i,j)}_t&=\left(\frac{(\tilde\chi^{(i)}_t-\tilde\chi^{(j)}_t)^2}{\tilde\chi^{(i)}_t\tilde\chi^{(j)}_t}\frac{\chi^{(i)}_t\chi^{(j)}_t}{(\chi^{(i)}_t-\chi^{(j)}_t)^2}
\right)^{\frac{\rho_i\rho_j}{4\kappa}}h'_t(0)^{\frac{\rho_i\rho_j}{4\kappa}}\\ 
\end{align*}
Note that all fractions are real numbers. Then
$(M^\varnothing\prod_iM^{(i)}\prod_{i<j}M^{(i,j)})$ is defined up to
some random positive time (possibly infinite), and is a local martingale.
\end{Lem}
\begin{proof}
This is a rather straightforward transposition to the radial case of results
and methods in \cite{LSW3}, which we discuss for the sake of
completeness. We use freely a complex-variable version of It\^o's formula.

Let $\tilde g_t=\phi_{\phi_A(K_t)}$, with the usual notations, so that
$h_t\circ g_t=\tilde g_t\circ \phi_A$. Then $(\tilde g_t)$ is a
time-changed radial Loewner chain:
$$\partial_t\tilde g_t=-\tilde g_t\frac{\tilde g_t+\tilde\xi_t}{\tilde
g_t-\tilde\xi_t}da_t$$
and $da_t=(h'_t(\xi_t)\xi_t/\tilde\xi_t)^2$. This is the key
``commutative diagram'' argument of \cite{LSW3}. Then, standard
differential calculus yields:
\begin{align*}
\partial_t h_t&=-h_t\frac{h_t+\tilde \xi_t}{h_t-\tilde
\xi_t}\partial_ta_t+h'_t.z\frac{z+\xi_t}{z-\xi_t}\\
[\partial_t
h_t](\xi_t)&=3\xi_th'_t(\xi_t)+3\xi_t^2h''_t(\xi_t)-3\frac{(\xi_th'_t(\xi_t))^2}{\tilde\xi_t}\\
\partial_th'_t&=-h'_t\frac{h_t+\tilde \xi_t}{h_t-\tilde
\xi_t}\partial_ta_t+h_th'_t\frac{2\tilde\xi_t}{(h_t-\tilde\xi_t)^2}\partial_ta_t+h_t''.z\frac{z+\xi_t}{z-\xi_t}+h'_t.\left(\frac{z+\xi_t}{z-\xi_t}-\frac{2z\xi_t}{(z-\xi_t)^2}\right)\\
[\partial_th'_t](\xi_t)&=\left(\frac 43\xi_t^2h_t'''-\frac{(\xi_th''_t)^2}{2h'_t}+3\xi_th''_t-\frac{\xi_t^2(h'_t)^3}{\tilde\xi_t^2}+h'_t\right)(\xi_t)\\
\partial_th'_t(0)&=h'_t(0)(\partial_ta_t-1)
\end{align*}
Relax the assumptions for the moment, and assume that $\xi$ satisfies
the SDE:
$$d\xi_t=i\xi_t\sqrt\kappa dB_t-\frac\kappa2\xi_tdt+\xi_tb_tdt$$
where $B$ is a standard Brownian motion and $b_t$ is a (progressive) drift coefficient.
Define:
\begin{align*}
\sigma_t&=\frac{\xi_th''_t(\xi_t)}{h'_t(\xi_t)}+1-\frac{\xi_th'_t(\xi_t)}{\tilde\xi_t}\\
\sigma^{(i)}_t&=\frac{\chi^{(i)}_t+\xi_t}{\chi^{(i)}_t-\xi_t}-\frac{\tilde\chi^{(i)}_t+\tilde\xi_t}{\tilde\chi^{(i)}_t-\tilde\xi_t}\frac{\xi_th'_t(\xi_t)}{\tilde\xi_t}
\end{align*}
Then:
\begin{align*}
\frac
{dM^\varnothing_t}{M^\varnothing_t}&=\frac{6-\kappa}{2\kappa}\sigma_t\left(i\sqrt\kappa
dB_t+b_tdt\right)\\
\left(\frac{\rho_i}{2\kappa}\right)^{-1}\frac{dM^{(i)}_t}{M^{(i)}_t}&=
\sigma^{(i)}_t\left(i\sqrt\kappa
dB_t+b_tdt+\frac{\rho_i}{2}\frac{\xi_t+\chi^{(i)}_t}{\xi_t-\chi^{(i)}_t}dt\right)-\frac{(6-\kappa)}{2}\sigma_t\frac{\tilde\chi^{(i)}_t+\tilde\xi_t}{\tilde\chi^{(i)}_t-\tilde\xi_t}\frac{\xi_th'_t(\xi_t)}{\tilde\xi_t}dt\\
\left(\frac{\rho_i\rho_j}{4\kappa}\right)^{-1}\frac{\partial_tM^{(i,j)}_t}{M^{(i,j)}_t}&=
\sigma^{(i)}_t\sigma^{(j)}_t+\frac{\xi_t+\chi^{(i)}_t}{\xi_t-\chi^{(i)}_t}\sigma^{(j)}_t+\frac{\xi_t+\chi^{(j)}_t}{\xi_t-\chi^{(j)}_t}\sigma^{(i)}_t
\end{align*}
Given that
$b_t=-\sum_i\frac{\rho_i}2\frac{\xi_t+\chi^{(i)}_t}{\xi_t-\chi^{(i)}_t}$,
if $N$ denotes $(M^\varnothing\prod_iM^{(i)}\prod_{i<j}M^{(i,j)})$, one gets:
\begin{align*}
\frac{dN_t}{N_t}&=\frac{dM^\varnothing_t}{M^\varnothing_t}+\sum_i\left(\frac{dM^{(i)}_t}{M^{(i)}_t}+\frac{d\langle
M^\varnothing_t,M^{(i)}_t\rangle}{M^\varnothing_tM^{(i)}_t}\right)
+\sum_{i<j}\left(\frac{dM^{(i,j)}_t}{M^{(i,j)}_t}+\frac{d\langle
M^{(i)}_t,M^{(j)}_t\rangle}{M^{(i)}_tM^{(j)}_t}\right)\\
&=\sigma'_tdB_t+\frac{6-\kappa}{2\kappa}b_t\sigma_t+\sum_i\left(\frac{6-\kappa}{2\kappa}\sigma_t\frac{\rho_i}2\frac{\xi_t+\chi^{(i)}_t}{\xi_t-\chi^{(i)}_t}+\frac{\rho_i}{2\kappa}\sigma^{(i)}_t\left(b_t+\frac{\rho_i}{2}\frac{\xi_t+\chi^{(i)}_t}{\xi_t-\chi^{(i)}_t}\right)\right)dt\\
&\hphantom{=}+\sum_{i<j}\frac{\rho_i\rho_j}{4\kappa}\left(\frac{\xi_t+\chi^{(i)}_t}{\xi_t-\chi^{(i)}_t}\sigma^{(j)}_t+\frac{\xi_t+\chi^{(j)}_t}{\xi_t-\chi^{(j)}_t}\sigma^{(i)}_t\right)dt\\
&=\sigma'_tdB_t
\end{align*}

\end{proof}

As in the chordal case, the symmetry of these formulae when
$\underline\rho=(2,\dots,2)$ enables to derive restriction formulae
for $n$-braid SLEs.

\begin{Lem}
Let $(K.)$ be a radial $n$-braid SLE in $\U$ started from distinct points
$\chi_1,\dots,\chi_n$, and $A$ be a hull not intersecting these
points. If $\kappa\leq 8/3$, and $L$ is an independent loop soup with
intensity $\lambda_\kappa$, then:
$$\P(K_\infty^L \cap
A=\varnothing)=|\phi_A'(0)|^{2h_{0;n/2}}\prod_i|\phi_A'(\chi_i)|^{h_{1;2}}\prod_{i<j}\left|\frac{\phi_A(\chi_j)-\phi_A(\chi_i)}{\chi_j-\chi_i}\right|^{2/\kappa}$$
\end{Lem}

If $\chi^{(1)}_0,\dots,\chi^{(n)}_0$ collapse to $\chi$, the above formula reduces to:
$$|\phi_A'(\chi)|^{h_{1;n+1}}|\phi_A'(0)|^{2h_{0;n/2}}$$

\begin{Cor}[Restriction property]
Let $\chi_1,\dots,\chi_n$ and $K_\infty^L$ as above. This defines a family of
probability measures $\mu_{\kappa,\underline\chi}$ on ``sea star''
hulls $K'=K^L$. %(if $n=5$; ``Eiffel tower'' hulls if $n=3$ \dots). 
This family has the
restriction property:\\
for all hull $A$, $(\phi_A)_*\mu_{\kappa,\underline\chi}(.|K'\cap A=\varnothing)=\mu_{\kappa,(\phi_A)_*(\underline\chi)}$
\end{Cor}
\begin{proof}
As in the chordal case.
\end{proof}

\section{Multiply connected domains}

One can think of SLE as a diffusion in a configuration space. The
diffusion coefficients are constrained by the conformal invariance
requirement. If the associated moduli space is a point, then the
coefficients are constant parameters; this situation corresponds to
chordal and radial SLE, and $\SLE(\kappa,\rho)$ (and also ``annulus
SLE''). If the moduli space
is larger, then SLE is essentially specified by the data of diffusion
coefficients as functions on the moduli space; so we are no longer in
the parametric situation. We only discuss the ``constant $\kappa$''
case, for physical and technical reasons. Also, we will be mainly
interested in expressing necessary conditions for reversibility in
multiply connected domains, so we will not carry the discussion in the
same degree of generality as in the simply connected case.

In the case of a simply connected domains with $2n$ points marked on
the boundary, $\SLE$ is specified by $\kappa$ and a function of $(2n-3)$
independent cross-ratios of the $2n$ boundary points. If we add the
requirement that the $\SLE$ commutes with $(2n-1)$ $\SLE$s started at the
other points, then we have to choose a ``partition function'' $\psi$
as discussed earlier. This function belongs to the finite-dimensional solution
space of a holonomic system derived from the commutation conditions,
and the situation is parametric again ($1+C(2n,n)$ parameters).

Similarly, in the case of multiply connected domain, we want to
restrict the diffusion coefficients to the ``physically relevant''
ones.
We consider in particular the case of chordal SLE (going from $x$ to
$y$, $x$ and $y$ on the same boundary component) in a multiply
connected domain.

There are at least two ways to describe SLE in a multiply connected
domain. The first one, that follows closely the simply connected case,
consists in choosing a parametric family of standard domains (a section of the
moduli space), and writing explicit diffusion equations for the
parameters; this is the approach of \cite{D3,BauFr}. Another route,
following Makarov and Zhan (see \cite{DZT}), consists in using a local chart at
the growth point and a ``conformally invariant SDE'', so that the path
distribution does not depend on the choice of local chart. In the first case,
the diffusion coefficient is a function on the moduli space; in the
other case, SLE is specified by a ``partition function'', which is a
conformally covariant function on the configuration space; taking its
log derivative (w.r.t. the  growth point), one gets a function on the
moduli space. We will use this second framework, that better suits our
purposes.

So let ${\mf C}$ be the configuration space of $(g+1)$-connected
plane domains with $(m+2)$ points marked on the boundary and $n$ points
marked in the bulk. Two of the marked points, $x$ and $y$ are on the same
component of the boundary. Denote by ${\mf M}$ the associated moduli
space.

First we briefly summarize the local chart approach. Any configuration
is equivalent to a configuration of type $(\H\setminus K,x,\dots)$ where $x$ is
real and $K$ is a compact subset of $\H$ (with $g$ connected component). By conformal invariance, we need
only to define SLE for these configurations, and need to do it
coherently (independently of choices). Let $h$ be a conformal
equivalence between $c=(\H\setminus K,x,y=\infty,\dots)$ and
$h_*c=(\H\setminus K',x',y'=\infty,\dots)$. SLE in
$c$ is defined by the chordal Loewner equations and an SDE:
$$\partial_t g_t=\frac 2{g_t-W_t}, dW_t=\sqrt\kappa dB_t+b(c_t)dt$$
From \cite{LSW3}, we can write the SDE for the driving process of
the image of the SLE by $h_0=h$:
$$dh_t(W_t)=h'_t(W_t)dW_t+\left(\frac\kappa 2-3\right)h''_t(W_t)dt$$
where $h_t=\tilde g_t\circ h\circ g_t^{-1}$, and $\tilde g_t$ also
solves the chordal Loewner equations (though with a time
change). After a time change, and at time 0, one sees that the
condition necessary for invariance of the SDE is the following
covariance condition:
$$b(h_*c)=\frac{b(c)}{h'(x)}+\left(\frac\kappa 2-3\right)\frac{h''}{h'^2}(x).$$
Here $h$ is normalized by $h(\infty)=\infty$, $h'(\infty)=1$
(hydrodynamic normalization at infinity).
Let $\psi$ be a positive function on the configuration space. We say
that $\psi$ is $\alpha$-covariant if:
\begin{enumerate}
\item (M\"obius invariance) For any $c=(D\setminus K,\dots)$, $c'=(D'\setminus K',\dots)$
where $D,D'$ are simply connected, $h:c\rightarrow c'=h_*c$ is an
equivalence of configurations, and $h$ extends to a conformal
equivalence $D\rightarrow D'$, one has:
$$\psi(h_*c)=\psi(c)$$
\item (covariance)  For any $c=(D\setminus K,x,y,\dots)$, $c'=(D\setminus K',x,y,\dots)$
where $D$ is simply connected, $x,y\in\partial D$, $\partial D$ is
smooth at $x,y$, and $h:c\rightarrow c'=h_*c$ is an
equivalence of configurations, one has:
$$\psi(h_*c)=(h'(x)h'(y))^{-\alpha}\psi(c)$$
\end{enumerate} 

Note that the function
$\psi$ is completely determined by these conditions and its restriction
to a section of the moduli space. Let us give three (important) examples
of such covariant functions, say for annuli with two marked points on one
component of the boundary: $c=(D,x,y)$. 
\begin{enumerate}
\item
In $c=(D,x,y)$, assume that the arc $(xy)$
is blue and $(yx)$ is yellow. Let $\psi(c)$ be the probability that
$(xy)$ and $(yx)$ are connected to the other boundary component by a
blue (resp. yellow) cluster in the scaling limit of critical
percolation. Alternatively, $\psi(c)$ is the corresponding $\SLE_6$
probability (see \cite{D3}). Then $\psi$ is $0$-covariant. (Also,
$\psi=1$ is $0$-covariant; this is a version of locality for $\SLE_6$).
\item 
In his thesis, Beffara uses the results of \cite{LSW3} and an
inclusion-exclusion argument to prove the following: let $\psi(c)$ be
the probability that chordal $\SLE_{8/3}$ from $x$ to $y$ in the
(filled) domain avoids the hole (resp. leaves the hole on its left, 
resp. leaves the hole on its right). Then $\psi$ is $5/8$-covariant. 
\item 
If $\kappa=2$, $\alpha_\kappa=1$, and $\SLE_2$ is the scaling
limit of Loop-Erased Random Walks (\cite{LSW2,DZ}); these walks are
closely related to some discrete harmonic quantities. Define:
$$\psi(c)=\frac{\partial^2}{\partial n_x\partial n_y} G_D(w,w')$$
the normal derivative at $x$ and $y$ of the Green kernel (which is
symmetric in the two variables). Then the invariance property of the
Green kernel implies that $\psi$ is $1$-covariant. In general domains,
this is the (chordal version of) Harmonic Random Loewner Chain (HRLC)
as defined by
Zhan in \cite{DZT}. Similar harmonic constructions exist for
$\kappa=8$, $\alpha_\kappa=-1/8$, using normal reflection on some
boundary components.
\end{enumerate}

If $c=(\H\setminus K,x,y=\infty,\dots)$, $h_*c=(\H\setminus
K',x',y'=\infty,\dots)$ are equivalent configurations, $h'(\infty)=1$,
and
we define
$b=\kappa\partial_x\psi/\psi$, we get:
$$b(h_*c)=\frac\kappa{h'(x)}\frac{\partial_x(h'(x)^{-\alpha}\psi(c))}{h'(x)^{-\alpha}\psi(c)}=\frac{b(c)}{h'(x)}-\kappa\alpha\frac{h''}{h'^2}(x)$$
where $c$ is implicitly a function of $x$ (everything else being fixed).
So if $\alpha=\alpha_\kappa=h_{1;2}(\kappa)=(6-\kappa)/2\kappa$, one
can define an SLE starting from the $\alpha$-covariant partition function
$\psi$ (at least up to some positive stopping time). We denote this
by $\SLE_\kappa(\psi)$. This is well-defined for some positive time; we will not consider here the (difficult)
questions of long-time behaviour.

\subsection{Commutation conditions}

Now assume we are given two $\alpha_\kappa$-covariant functions
$\psi_1$, $\psi_2$ on the configuration space, and use
them to define SLEs starting resp. at $x$ and $y$. 
%(once again another arbitrary point on the same boundary component is marked for normalization
purposes). 
In a configuration $c=(\H\setminus K,x,y,\dots)$, one grows $\SLE$s at
$x$ and $y$ up to capacity $\eps$, $\eps'$ (seen from infinity in
$\H$; this is also arbitrary) and consider the effect on functions of
$x$ and $y$ (after erasing hulls using chordal SLE). As earlier, this
leads to the commutation conditions on $\psi_1$,
$\psi_2$. To make the argument neater, we compute on a section of the
moduli space, as in \cite{D3,BauFr}.  As before, there is a marked
point on the boundary used for normalization.
In fact, in this chordal setup,
it is more convenient not to quotient by automorphisms that fix this point; so all
the construction will commute with scaling and translation.
For definiteness, consider the following type of configurations: the
upper half-plane $\H$ minus horizontal slits with appropriate marked
points (including $\infty$). The only equivalences between such
configurations are given by scaling and translation. Let ${\mf D}$ be
this family. 

Consider an element $D_0=(H_0,x,y,\dots)\in {\mf D}$, that is $H_0$ is
the half-plane minus some horizontal slits. We grow an
$\SLE_\kappa(\psi_1)$ at $x$ up to half-plane capacity $\eps$, then an 
$\SLE_{\kappa}(\psi_2)$ at $y$ (in the remaining domain) stopped when it
reaches capacity $\eps$ (in the original domain). We then revert the
order of the 
procedure and look for necessary conditions for these two procedures
to yield the same distribution (of configuration). In particular, we
compare the effect on moduli, up to second order in $\eps$.

So consider a Loewner chain $(K_t)$ growing at $x$. Let $g_t$ be a conformal
equivalence $\H\setminus K_t\rightarrow \H$, and $f_t$ a conformal
equivalence $D_0\setminus K_t\rightarrow D_t$ for some $D_t\in{\mf
D}$. Everything is uniquely defined if we impose hydrodynamic
normalization at infinity for $g_t,f_t$. Then let $h_t=f_t\circ g_t^{-1}$. By
construction, $g_t$ solves the Loewner equations:
$$\partial_t g_t=\frac 2{g_t-X_t}, dX_t=\sqrt\kappa
dB_t+\kappa\frac{\partial_x\psi_1(H_t)}{\psi_1(H_t)}dt$$
where $\psi_1$ is evaluated at the configuration $H_t=(g_t(H_0\setminus
K_t),X_t,\dots)$. Now consider $(\partial_t f_t)\circ f_t^{-1}$. This is a meromorphic
function on $D_t$, taking real values on $\R$, vanishing at $\infty$,
regular except at $x$ where it has a simple pole, and with constant
imaginary part on the slits. It is easy to see that for each $D\in
{\mf D}$, there is a unique function $V_D$ (Schwarz kernel) satisfying these conditions
and with residue $2$ at $x$. Also define $A_D,B_D$ by:
$$V_D(w)=\frac 2{w-x}+A_D+B_D(w-x)+O((w-x)^2)$$ 
Then:
\begin{align*}
\partial_t f_t&=h'_t(X_t)^2 V_{D_t}\circ f_t\\ 
\partial_t h_t(w)&=\partial_t (f_t\circ g_t^{-1})(w)=h'_t(X_t)^2 V_{D_t}\circ h_t(w)-\frac{2h'_t(w)}{w-X_t}\\
\partial_t h'_t(w)&=h'_t(X_t)^2.(h'_t.V'_{D_t}\circ h_t)(w)-\frac{2h''_t(w)}{w-X_t}+\frac{2h'_t(w)}{(w-X_t)^2}.
\end{align*}
Taking limits at $w=X_t$, one gets:
\begin{align*}
(\partial_t h_t)(X_t)&=h'_t(X_t)^2A_{D_t}-3h''_t(X_t)\\
(\partial_t h'_t)(X_t)&=h'_t(X_t)^3B_{D_t}+\left(\frac {h''_t}{2h'_t}-\frac{4h'''_t}3\right)(X_t)
\end{align*}
which gives in particular
\begin{align*}
dh_t(X_t)&=h'_t(X_t)dX_t+\left(\frac\kappa
2-3\right)h''_t(X_t)dt+h'_t(X_t)^2A_{D_t}dt\\
dh'_t(X_t)&=h''_t(X_t)dX_t+\left(\frac
{h''_t(X_t)^2}{2h'_t(X_t)}+\left(\frac\kappa 2-\frac 43\right)h'''_t(X_t)\right)dt+h'_t(X_t)^3B_{D_t}dt
\end{align*}
as in \cite{LSW3} (and \cite{D3} for the case of annuli, using the
explicit Villat kernel). Also:
$$\partial_t f'_t(y)=h'_t(X_t)^2. (f'_t. V'_{D_t}\circ f_t)(y).$$
Now ${\mf
D}$ can be parametrized by a list of complex numbers
$(x,y,z_1,\dots,z_m)$ containing marked points on $\R$, other marked
points (including endpoints of horizontal slits) and their
conjugates. Then the infinitesimal generator at $D$ is:
$${\mc L}_1=\frac\kappa 2\partial_{xx}
+\left(A_D+\kappa\frac{\partial_x\psi_1}{\psi_1}\right)\partial_x
+V_D(y)\partial_y+\sum_i V_D(z_i)\partial_{z_i}$$
with $V_D=V_{(x,y,z_1,\dots)}$ and $\psi_1$ is evaluated at $D$; the
restriction
of $\psi_1$ to ${\mf D}$ can also be seen as a function of $(x,y,z_1,\dots)$.
The generator ${\mc L}_2$ associated with the $\SLE_{\kappa}(\psi_2)$
growing at $y$ is derived in similar fashion, and we get the
commutation condition:
$$[{\mc L}_1,{\mc L}_2]=2V'_y(x){\mc L_1}-2V'_x(y){\mc L_2}$$
where $V_x(z)=2/(z-x)+\cdots$ and $V_y(z)=2/(z-y)+\cdots$ are the
Schwarz kernels with poles at $x$ and $y$ respectively (and depend
implicitly on the other moduli). We will use the notation:
$$\ell_x=V_x(y)\partial_y+\sum V_x(z_i)\partial_i,{\rm\ \ } 
\ell_y=V_y(x)\partial_x+\sum V_y(z_i)\partial_i$$ 
for the part that does not involve $\kappa$, $\psi_1$, $\psi_2$.

Expanding the commutation condition
, we get the equations:
\begin{equation}\label{cmu1}
\partial_{xy}\log(\psi_1/\psi_2)=0
\end{equation}
and:
\begin{equation}\label{cmu2}
\begin{split}
-{\mc L_2}\left(A_x+\kappa\frac{\partial_x\psi_1}{\psi_1}\right)+{\mc
L}_1(V_y(x))+2V'_x(y)V_y(x)-2V'_y(x)\left(A_x+\kappa\frac{\partial_x\psi_1}{\psi_1}\right)&=0\\
-{\mc L_1}\left(A_y+\kappa\frac{\partial_y\psi_2}{\psi_2}\right)+{\mc
L}_2(V_x(y))+2V'_y(x)V_x(y)-2V'_x(y)\left(A_y+\kappa\frac{\partial_y\psi_2}{\psi_2}\right)&=0
\end{split}\end{equation}
All other
conditions are identities not involving $\psi_1$, $\psi_2$:
\begin{equation}\label{cmu3}
(\ell_x+2V'_x(y))V_y(z)=(\ell_y+2V'_y(x))V_x(z)
\end{equation}
For this, note that $V_y(z)$ does not depend on $x$ and other marked
points, but depends on endpoints of horizontal slits. Considering the
difference between left-hand side and right-hand side as a function of $z$, in particular its
expansion at $x$, $y$, one sees that it extends to a
bounded holomorphic function on the Schottky double of $D$ and
vanishes at infinity, hence is identically 0. The probabilistic
interpretation of the Schwarz kernel $V$ given in \cite{Law2} can also
be used to prove this identity, along the lines of Lemma 14 (ii) below.

Consider the conditions (\ref{cmu1}),(\ref{cmu2}). 
From (\ref{cmu1}), one can write $\psi_1=\psi_2=\psi$, by multiplying
$\psi_1$ by a function that does not depend on $x$ (so that the
definition of $\SLE_\kappa(\psi_1)$ is not affected), and doing the
same for $\psi_2$. One can write:
\begin{align*}
{\mc
L}_1\left(\frac{\partial_y\psi}{\psi}\right)-\partial_y\left(\frac{{\mc
L}_1\psi}{\psi}\right)
&=\frac{[{\mc L}_1,\partial_y]\psi}{\psi}+\kappa\frac{\partial_y\psi(\partial_x\psi)^2}{\psi^3}-\kappa\frac{\partial_x\psi\partial_{xy}\psi}{\psi^2}\\
&=-V'_x(y)\partial_y\log(\psi)-\kappa\partial_{xy}\log(\psi)\partial_x\log(\psi)+\kappa\partial_y\log(\psi)(\partial_x\log(\psi))^2-\kappa\frac{\partial_x\psi\partial_{xy}\psi}{\psi^2}\\
&=-V'_x(y)\partial_y\log(\psi)-\kappa\partial_y\left(\partial_x\log\psi\right)^2
\end{align*}
Let $\check{\mc L}_1$ be the differential operator obtained by setting
$\psi_1=1$ in ${\mc L}_1$:
$$\check{\mc L}_1=\frac\kappa 2\partial_{xx}
+A_x\partial_x
+V_D(y)\partial_y+\sum_i V_D(z_i)\partial_{z_i}=\frac\kappa 2\partial_{xx}
+A_x\partial_x
+\ell_x$$
and $\tilde {\mc L_2}$ is defined in the same fashion. Then (\ref{cmu2})
can be written as:
$$\kappa\partial_y\left(\frac{\check{\mc L}_1\psi}\psi\right)=\check{\mc L}_2(V_x(y))-\check{\mc L}_1(A_y)+2V'_y(x)V_x(y)-2A_yV'_x(y)$$
and similarly for the other condition. Consider now (\ref{cmu3}),
applied to a marked point $z=y+\eps$:
\begin{eqnarray*}
0&=&V_x(y)\partial_yV_y(z)+V_x(z)\partial_zV_y(z)+\sum_iV_x(z_i)\partial_{z_i}V_y(z)-V_y(x)\partial_xV_x(z)-V_y(z)\partial_zV_x(z)-\sum_iV_y(z_i)\partial_{z_i}V_x(z)\\
&&+2V'_x(y)V_y(z)-2V'_y(x)V_x(z)\\
&=&V_x(y)\left(\frac
2{\eps^2}+\partial_yA_y-B_y+(\partial_yB_y-2C_y)\eps\right)+\left(V_x(y)+\eps
V'_x(y)+\frac{\eps^2}2V''_x(y)+\frac{\eps^3}6V'''_x(y)\right)\left(-\frac
2{\eps^2}+B_y+2C_y\eps\right)\\
&&+\sum_iV_x(z_i)\partial_{z_i}(A_y+\eps B_y)
-V_y(x)\partial_xV_x(y)-V_y(x)\partial_xV'_x(y)\eps-\left(\frac
2\eps+A_y+B_y\eps\right)\left(V'_x(y)+\eps V''_x(y)+\frac{\eps^2}2
V'''_x(y)\right)\\
&&-\sum_iV_y(z_i)\partial_{z_i}(V_x(y)+\eps V'_x(y))
+2V'_x(y)\left(\frac
2\eps+A_y+B_y\eps\right)-2V'_y(x)V_x(y)-2V'_y(x)V'_x(y)\eps+O(\eps^2)
\end{eqnarray*}
where $C_y$ is such that
$V_y(z)=2/\eps+A_y+B_y\eps+C_y\eps^2+O(\eps^3)$.
Considering the coefficients of $\eps^0$, $\eps^1$, it follows that:
%\begin{align*}
%0&=V_x(y)\partial_y A_y-3V''_x(y)+\sum_iV_x(z_i)\partial_{z_i}A_y-V_y(x)\partial_xV_x(y)+V'_x(y)A_y
%-\sum_iV_y(z_i)\partial_{z_i}V_x(y)-2V'_y(x)V_x(y)\\
%&=\check{\mc L}_1(A_y)-\check{\mc L}_2(V_x(y))+\left(\frac{\kappa} 2-3\right)V_x''(y)+2V'_x(y)A_y-2V_x(y)V'_y(x)
%\end{align*}
\begin{equation}\label{cmu4}
(\ell_x+2V'_x(y))(A_y)-(\check{\mc
L}_2+2V'_y(x))(V_x(y))+\left(\frac{\kappa} 2-3\right)V_x''(y)
=0
\end{equation}
and
\begin{equation}\label{cmu5}
(\ell_x+2V'_x(y))(B_y)-(\check{\mc
L}_2+2V'_y(x))(V'_x(y))+\left(\frac{\kappa} 2-\frac 43\right)V_x'''(y)
=0
\end{equation}
Of course one can exchange the roles of $x$ and $y$ in these
identities. So (\ref{cmu2}) can be written as:
$$\kappa\partial_y\left(\frac{\check{\mc L}_1\psi}\psi\right)=
\left(\frac{\kappa} 2-3\right)V_x''(y)$$
or (with the symmetric condition):
$$\partial_y\left(\frac{\left(\check{\mc L}_1+\alpha_{\kappa} V'_x(y)\right)\psi}\psi\right)=\partial_x\left(\frac{\left(\check{\mc L}_2+\alpha_\kappa V'_y(x)\right)\psi}\psi\right)=0.$$

\subsection{Restriction martingales}

Restriction-like martingales in multiply-connected domains involving
harmonic invariants are studied in
\cite{Law2,DZT}, in particular when $\kappa=2$. When $\kappa\neq 2$,
these ``harmonic'' martingales are distinct from those we are
discussing in these sections, though methods are fairly similar.

Let $D$ be a subdomain of $\H$ (which we can freely assume if we want
to define M\"obius invariant distributions); the boundary of $D$
contains open real segments around $x,y\in\R$.
Define: $\Gamma(D,x)=\mu^{bub}_x(\{\delta\subsetneq D\})$, where
$\mu^{bub}_x$ is the bubble measure rooted at $x$ (see \cite{LW}). Then:
$$\Gamma(D,x)=\partial_{n_y}\left(M_\H(y,x)- M_D(y,x))\right)_{|y=x}$$
%$$\Gamma(D,x)=\pi\partial_{n_y}\left(\partial_{n_x}(G_\H(y,x)- G_D(y,x))\right)_{|y=x}$$
where $M_D$ denotes the minimal function with singularity at $x$,
i.e. the positive harmonic function in $D$ that extends continuously to 0 on
the boundary except at $x$, with normalization
$M_D(x+i\eps,x)=\eps^{-1}(1+O(\eps))$. This function can be obtained by
taking the normal derivative at $x$ of the Green's function $G_D(y,x)$
(with adequate normalization). 
Let $\phi$ be a
conformal equivalence $D\rightarrow D'=\phi(D)$.
From the conformal
invariance property of the Green's function, it is easy to derive the
covariance property of the minimal function:
$$M_D(y,x)=\phi'(x)M_{\phi(D)}(\phi(y),\phi(x))$$
and then the Schwarzian-like covariance property for $\Gamma$: 
$$\Gamma(D,x)=\phi'(x)^2\Gamma(\phi(D),\phi(x))-\frac{S\phi(x)}6.$$
Note that in the case where $D$ is simply connected,
$\Gamma(D,x)=-S\phi_D(x)/6$, and this is the usual covariance property
of Schwarzian derivatives: 
$$S(\phi_D)=S(\phi_{D'}\circ\phi)=(\phi')^2(S\phi_{D'})\circ\phi+S\phi$$

We are interested in the following situation: let $\gamma$ be a
chordal $\SLE_\kappa$ in $\H$, say $\kappa\leq 8/3$, and $L$ is an
independent loop soup with intensity $\lambda_\kappa$. Condition on
the event that $\gamma^L=\gamma\cup\{\delta\in
L:\delta\cap\gamma\neq\varnothing\}$ stays in $D_0$. It is not clear
at this point that the probability of this event is $\alpha_\kappa$-
covariant, hence that the resulting distribution is conformally
invariant.

First consider a chordal $\SLE_\kappa$ from $x$ to $y$ in $\H\supset
D_0$ (a chordal $\SLE$ unaware of the presence of holes), and an
$\alpha_\kappa$-covariant function $\varphi$ on the moduli space. By M\"obius
invariance, one can send $y$ to infinity (and use hydrodynamic
normalization); so $(X_t/\sqrt\kappa)$ is now a standard Brownian
motion. Here $g_t:D_0\setminus K_t\rightarrow H_t$ is a conformal
equivalence that extends through the holes, $h_t:H_t\rightarrow
D_t$ is a conformal equivalence and $D_t\in{\mf D}$ is a standard domain.
Let $M_t=\varphi(D_0\setminus K_t,\dots)$. From the covariance assumption,
$$M_t=h'_t(X_t)^\alpha\varphi(h_t(X_t),\dots).$$
Set $\alpha=\alpha_\kappa=(6-\kappa)/2\kappa$. Then:
$$\frac{dh'_t(X_t)^\alpha}{h'_t(X_t)^\alpha}=\alpha\frac{h''_t(X_t)}{h'_t(X_t)}dX_t-\frac{\lambda_\kappa}6
Sh_t(X_t)dt+\alpha h'_t(X_t)^2B_{D_t}dt$$
where $\lambda_\kappa=(8-3\kappa)(6-\kappa)/2\kappa$ and $Sh_t$ is the
Schwarzian derivative. Besides:
\begin{align*}
d\varphi(h_t(X_t),\dots,f_t(z_i),\dots)=&h'_t(X_t)\partial_x\varphi dX_t+\left[\frac\kappa 2 h'_t(X_t)^2\partial_{x,x}+
\left(\frac\kappa
2-3\right)h''_t(X_t)\partial_x\right.\\
&\left.+h'_t(X_t)^2A_{D_t}\partial_x+h'_t(X_t)^2\sum_iV_{D_t}(z_i)\partial_{z_i}\right]\varphi(h_t(X_t),\dots,f_t(z_i),\dots)dt
\end{align*}

Let:
$$Z_t=M_t\exp\left(-\lambda_\kappa\int_0^t\Gamma(H_s,X_s)ds\right)=M_t\exp\left(-\lambda_\kappa\int_0^t\left(h'_s(X_s)^2\Gamma(D_s,h_s(X_s))-\frac
{Sh_s(X_s)}6\right)ds\right)$$ 
Then $Z_t$ is a local martingale iff $\varphi$
(restricted to ${\mf D}$, hence considered as a function of the
parameters $x,\dots,z_i,\dots$) is annihilated by the differential operator:
$${\mc
M}=\frac\kappa 2\partial_{xx}+A_D\partial_x+\ell_x+\alpha
B_D-\lambda_\kappa\Gamma_x$$
where $\Gamma_x=\Gamma(D,x)$ in the standard domain $D$. 
If $y$ is finite (and another marked point is used for normalization),
one can compute along the same lines. So if $C_t$ is the covariance
factor:
$$C_t=\left(h'_t(X_t)h'_t(Y_t)\left(\frac{Y_t-X_t}{h_t(Y_t)-h_t(X_t)}\right)^2\right)^\alpha$$
then
\begin{eqnarray*}
\frac {dC_t}{C_t}&=&\left(\alpha\frac{h''_t(X_t)}{h'_t(X_t)}dX_t-\frac{\lambda_\kappa}6
Sh_t(X_t)dt+\alpha h'_t(X_t)^2B_{D_t}dt\right)\\
&&+\alpha\left(h'_t(X_t)^2V'_{D_t}(h_t(Y_t))+\frac{2}{(Y_t-X_t)^2}\right)dt
+2\alpha\left(\frac{dX_t}{X_t-Y_t}+(5-\kappa)\frac{dt}{(Y_t-X_t)^2}\right)\\
&&-\frac{2\alpha}{h_t(X_t)-h_t(Y_t)}\left(h'_t(X_t)dX_t+\left(\frac\kappa
2-3\right)h''_t(X_t)dt+h'_t(X_t)^2(A_{D_t}-V_{D_t}(h_t(Y_t))-\frac{3}{h_t(X_t)-h_t(Y_t)})dt\right)\\
&&+2\kappa\alpha^2\left(\frac{h''_t(X_t)}{h'_t(X_t)(X_t-Y_t)}-\frac{h''_t(X_t)}{h_t(X_t)-h_t(Y_t)}-2\frac{h'_t(X_t)}{(X_t-Y_t)(h_t(X_t)-h_t(Y_t))}\right)dt\\
&=&\alpha\left(\frac{h''_t(X_t)}{h'_t(X_t)}+\frac{2}{X_t-Y_t}-\frac{2h'_t(X_t)}{h_t(X_t)-h_t(Y_t)}\right)\sqrt\kappa
dB_t-\frac{\lambda_\kappa}6
Sh_t(X_t)dt\\
&&+\alpha h'_t(X_t)^2\left(B_{D_t}+V'_{D_t}(h_t(Y_t))-\frac{2}{h_t(X_t)-h_t(Y_t)}(A_{D_t}-V_{D_t}(h_t(Y_t))-\frac{3}{h_t(X_t)-h_t(Y_t)})\right)dt
\end{eqnarray*}
So $M_t$ is a local martingale iff $\varphi$ is annihilated by the operator:
$${\mc
M}=\frac\kappa 2\partial_{xx}+(A_D+\frac{\kappa-6}{x-y})\partial_x+\ell_x+\alpha
\left(B_D+V'_D(y)+\frac{2(V_D(y)-A_D)}{x-y}+\frac{6}{(x-y)^2}\right)-\lambda_\kappa\Gamma_x$$
Observe that:
\begin{align*}
(x-y)^{-2\alpha}{\mc M}(x-y)^{2\alpha}&=\frac\kappa 2\partial_{xx}+A_D\partial_x+\ell_x+\alpha
\left(B_D+V'_D(y)\right)+\lambda_\kappa\Gamma_x\\
&=\tilde {\mc
L_1}+\alpha (B_D+V'_D(y))-\lambda_\kappa\Gamma_x
\end{align*}
The conjugation corresponds to a change of reference measure ($\SLE$
aiming at $\infty$ rather than $\SLE$ aiming at $y$).

We sum up the discussion of this section. As before, ${\mf C}$ is a
configuration space, $x$ and $y$ are marked points on a boundary
component of a configuration; ${\mf M}$ is the associated moduli
space, and ${\mf D}$ is a class of standard domains, with associated
Schwarz kernels $V$. We note $V_x$, $V_y$ to emphasize the pole of
the kernel.

\begin{Prop}
(i) Let $\psi_1$, $\psi_2$ be $\alpha_\kappa$-covariant functions;
consider an $\SLE_\kappa(\psi_1)$ growing at $x$ and an $\SLE_\kappa(\psi_2)$
growing at $y$. These two SLEs satisfy commutation condition iff there
is an $\alpha_\kappa$-covariant function $\psi$ such that
$\SLE_\kappa(\psi_i)$, $i=1,2$, is distributed as an $\SLE_\kappa(\psi)$ growing at
$x$ (resp. $y$) and $\psi$ satisfy the conditions:
$$\partial_y\left(\frac{\left(\check{\mc L}_1+\alpha V'_x(y)\right)\psi}\psi\right)=\partial_x\left(\frac{\left(\check{\mc L}_2+\alpha V'_y(x)\right)\psi}\psi\right)=0$$
(ii) Let $\varphi$ be an $\alpha_\kappa$-covariant function. Let $c\in
{\mf C}$ be a configuration. One can assume that $c=D_0$ is a
standard domain. Consider a chordal $\SLE_\kappa$ from $x$
to $y$ in $c$ with holes erased; let $(c_t)_t$ be the corresponding
family of configurations. Let $h_t$ be the map from $c_t$ to a
standard domain $D_t\in {\mf D}$, with hydrodynamic normalization. Let:
$$Z_t=\varphi(c_t)\exp\left(-\lambda_\kappa\int_0^t
\Gamma(c_s,\gamma_s)ds\right)$$
Then $Z$ is a local martingale iff $\varphi$ restricted to ${\mf D}$
is annihilated by the operator:
$${\mc M}=
(x-y)^{2\alpha}\left(\check {\mc L_1}+\alpha (B_x+V'_x(y))-\lambda_\kappa\Gamma_x\right)(x-y)^{-2\alpha}$$
(iii) In the situation of $(ii)$, consider the Girsanov transform of
chordal $\SLE$ by $Z$. Then the resulting process is an
$\SLE_\kappa(\psi)$, where $\psi$ is an $\alpha_\kappa$-covariant function
whose restriction to ${\mf D}$ is annihilated by the operator:
$$\tilde {\mc L_1}+\alpha_\kappa (B_x+V'_x(y))-\lambda_\kappa\Gamma_x$$
In particular, $\SLE_\kappa(\psi)$ started at $x$ and
$\SLE_\kappa(\psi)$ started at $y$ satisfy the commutation
conditions. 
\end{Prop}

\subsection{The case $\kappa=0$}

In this section, we address the following (non trivial) question: how
to define a (deterministic) reversible $\SLE_0$ in a multiply connected domain ?
This requires implicitly the domain Markov property and conformal
invariance. While there are many ways to construct conformally
invariant chords (e.g. as level lines or flow lines of harmonic
invariants), it is not so obvious to satisfy all conditions
simultaneously. We propose here a construction based on the
restriction ideas (and do not claim to prove anything rigorous about
it).

Consider a nice subdomain $H_0$ of $\H$ (i.e. $\overline{\H\setminus
H_0}$ is compact and has finitely many components), $x_0\in\R$,
$(x_0-\eps,x_0+\eps)\subset \partial H_0$ for some $\eps>0$; we attempt to
describe an $\SLE_0(x_0\rightarrow\infty)$ in $H_0$. The idea is to take
$\gamma$ a chordal $\SLE_\kappa(x_0\rightarrow\infty)$ in $\H$ conditionally on
$\gamma^L$ staying in $H_0$, where $L$ is an independent loop soup
with intensity $\lambda_\kappa$, and take the limit as $\kappa\searrow
0$. On the one hand, the unconditional chordal $\SLE_\kappa$ converges to a hyperbolic
geodesic in $\H$, which is not conformally invariant; on the other
hand the intensity of the loop soup diverges. We shall define a total
cost summing the large deviation rate for Brownian motion and the loop
soup term, which will give a conformally covariant functional.

To check conformal invariance, we will use comparisons of Loewner
chains in conformally equivalent domains as in \cite{LSW3} (see also
the previous subsections).
Let us fix some notations. Let $\phi_0: H_0\rightarrow\tilde H_0$ be a
conformal equivalence to $\tilde H_0$ (and $\tilde H_0$ satisfies the same properties
as $H_0$). Consider the Loewner flow:
$$\partial_t g_t=\frac 2{g_t-X_t}$$
that maps $H_0$ to $H_t$ (and extends to $\H$). This Loewner chain is
mapped by $\phi_0$ to a Loewner chain corresponding to the flow:
$$\partial_t \tilde g_t=\frac {2\phi_t'(X_t)^2}{\tilde g_t-\tilde X_t}$$
where $\phi_t=\tilde g_t\circ\phi_0\circ g_t^{-1}$ is a conformal
equivalence $H_t\rightarrow \tilde H_t$ (and $\tilde X_t=\phi_t(X_t)$ etc...). 

Consider now the functional on driving processes:
$$I(X)=I(X,H_0)=\frac 12\int_0^\infty(\partial_t X_t)^2dt+24\int_0^\infty
\Gamma(H_t,X_t)dt.$$
Informally, as $\kappa\searrow 0$, a driving process $X$ has weight
$\propto\exp(-I(X)/\kappa)$, the first term being the large deviations
rate for Brownian motion, the second term coming from the loop soup
conditioning ($\lambda_\kappa\sim 24/\kappa$). We want to prove that
$I(X)=I(\tilde X)+c$, where $\tilde X$ has time parameter:
$$s(t)=\int_0^t\phi_t'(X_u)^2du.$$
Taking into account the time change, we get:
$$I(\tilde X)=I(\tilde X,\tilde H_0)=\frac 12\int_0^\infty\left(\frac{\partial_t\tilde
X_t}{\phi'_t(X_t)^2}\right)^2\phi'_t(X_t)^2dt+24\int_0^\infty
\phi'_t(X_t)^2\Gamma(\tilde H_t,\tilde X_t)dt.$$
From the covariance of $\Gamma$, we have:
$$\Gamma(H_t,X_t)=\phi'_t(X_t)^2\Gamma(\tilde H_t,\tilde X_t)-\frac {(S\phi_t)(X_t)}6.$$
Also, comparisons of Loewner chains yield:
\begin{align*}
\partial_t\tilde X_t&=\phi'_t(X_t)\partial_tX_t-3\phi_t''(X_t)\\
\partial_t(\phi_t'(X_t))&=\phi''_t(X_t)\partial_tX_t+\left(\frac{\phi_t''}{2\phi_t'}-\frac
43\phi'''_t\right)(X_t)
\end{align*}
It follows that:
\begin{align*}
I(\tilde X)-I(X)&=\frac
12\int_0^\infty\left(-6\frac{\phi_t''}{\phi_t'}(X_t)\partial_tX_t+9\left(\frac{\phi_t''}{\phi_t'}\right)^2(X_t)\right)dt+4\int_0^\infty
S\phi_t(X_t)dt\\
&=\frac
12\int_0^\infty\left(-6\frac{\phi_t''}{\phi_t'}(X_t)\partial_tX_t-3\left(\frac{\phi_t''}{\phi_t'}\right)^2(X_t)+8\frac{\phi_t'''}{\phi_t'}(X_t)\right)dt=-3\int_0^\infty\frac{\partial_t(\phi_t'(X_t))}{\phi_t'(X_t)}dt\\
&=3\log(\phi_0'(x_0))
\end{align*}
since $\phi_t'(X_t)\rightarrow 1$ ($\phi_t$ has hydrodynamic
normalization at $\infty$). Thus we get the covariance relation:
$$I(\phi_*X,\phi_*H)=I(X,H)+3\log(\phi'(x_0)).$$
So provided that $I$ attains a minimum, the minimizing path is
conformally invariant and Markov. Also, one can minimize over a subset
of paths satisfying particular global topological conditions
(e.g. leaving a hole on the left).
Assuming that chordal $\SLE_\kappa$
in $\H$ is reversible, this is also reversible (the conditioning event
does not depend on the orientation of the paths). For a 
configuration $(H,x,\infty)$, denote:
$$\psi(H)=\sup_X\left(\exp^{-I(H,X)}\right)$$
which is $(-3)$-covariant. With the previous notations, one can write
the coherence condition:
$$\psi(H)=\sup_{X_{[0,t]}}\left(\psi(H_t)\exp\left(-\frac 12\int_0^t(\partial_s X_s)^2ds+24\int_0^t
\Gamma(H_s,X_s)ds\right)\right)$$
By conformal invariance, one can assume that $H$ is a standard
configuration. Let $h_t:H_t\rightarrow D_t$ be the normalized
conformal equivalence to a standard domain. Then:
$$\psi(D_t)=\psi(H_t)(h'_t(X_t))^{-3}$$
and $\partial_t(h_t(X_t))_{|t=0}=A_x$, $\partial_t(h'_t(X_t))_{|t=0}=B_x$.
Let $\overline X$ denote a maximizing driving function.
Expanding the
previous condition for small $t$, one gets:
$$(\ell_x+A_x\partial_x+(\partial_t\overline X_t)\partial_x+3B_x)\psi-\left(\frac
12(\partial_t \overline X_t)^2+24\Gamma_x\right)\psi=0$$
at time $t=0$ (where $\psi$ is seen as a function on standard configurations). Besides, the left-hand side is optimal for
$\partial_t\overline X_t$. This implies that:
$$\partial_t \overline X_t=\frac{\partial_x\psi}{\psi}{\rm\ \ and\ \ }
\frac{(\partial_x\psi)^2}{2\psi}+(\ell_x+A_x\partial_x+3B_x-24\Gamma_x)\psi=
0.
$$
One can also get this equation by considering:
$$\left(\frac\kappa 2\partial_{xx}+A_x\partial_x+\ell_x+\alpha_\kappa B_x-\lambda_\kappa\Gamma_x\right)\psi_\kappa=0$$
and then rewrite the equation for $(\psi_\kappa)^{\kappa}$ and take
the limiting equation as $\kappa\searrow 0$. Now if the target point
is not $\infty$ but a finite point $y$, one can proceed similarly. The
condition obtained on $\psi$ for growth at $x$, $y$ gives the
commutation condition.

\subsection{Towards a classification}

In simply connected domains, we obtained a complete classification of
commuting SLEs. In the multiply connected case, it appears to be
much more technical, so we shall only outline some elements.
First, if we don't assume {\em a priori} that the two SLEs have the
same $\kappa$, then it is not hard to check that the commutation
condition for an $\SLE_\kappa(\psi)$ growing at $x$ and an
$\SLE_{\tilde\kappa}(\psi)$ growing at $y$ writes:
$$\partial_y\left(\frac{\left(\check{\mc L}_1+\alpha_{\tilde\kappa} V'_x(y)\right)\psi}\psi\right)=\partial_x\left(\frac{\left(\check{\mc L}_2+\alpha_\kappa V'_y(x)\right)\psi}\psi\right)=0$$
where now $\check{\mc L}_2=\frac{\tilde\kappa}2\partial_{yy}+\cdots$. Also,
if $\tilde\kappa\neq\kappa$, the covariance condition for $\psi$ is
modified as follows:

\begin{enumerate}
\item (M\"obius invariance) For any $c=(D\setminus K,\dots)$, $c'=(D'\setminus K',\dots)$
where $D,D'$ are simply connected, $h:c\rightarrow c'=h_*c$ is an
equivalence of configurations, and $h$ extends to a conformal
equivalence $D\rightarrow D'$, one has:
$$\psi(h_*c)=\psi(c)$$
\item (covariance)  For any $c=(\H\setminus K,x,x',y,y',\dots)$, $c'=(\H\setminus K',\dots)$
where $D$ is simply connected, $x,y\in\partial D$, $\partial D$ is
smooth at $x,y$, and $h:c\rightarrow c'=h_*c$ is an
equivalence of configurations, one has:
$$\psi(h_*c)=\left(h'(x)h'(x')\left(\frac{x-x'}{h(x)-h(x')}\right)^2\right)^{-\alpha_\kappa}\left(h'(y)h'(y')\left(\frac{y-y'}{h(y)-h(y')}\right)^2\right)^{-\alpha_{\tilde\kappa}}\psi(c)$$
\end{enumerate} 
where $x',y'$ are new target points for the SLEs. 

We now revert to the case $\kappa=\tilde\kappa$ (and the two SLEs are
``aiming at each other'') as discussed earlier and further study the
commutation conditions. Consider the operators:
$${\mc M}_1=\check {\mc L_1}+\alpha V'_x(y)+h_1,{\mc M}_2=\check {\mc L_2}+\alpha V'_y(x)+h_2$$
where $\partial_y h_1=\partial_x h_2=0$. To define two commuting SLEs,
we have to find functions $h_1,h_2,\psi$ such that ${\mc M}_1\psi={\mc
M}_2\psi=0$.

Say $\kappa=8/3$, and consider a chordal $\SLE_{8/3}$ in a simply
connected domain conditioned to avoid some holes. Then the conditional
$\SLE$ can be represented as an $\SLE_{8/3}(\psi)$, and $\psi$
(restricted to a section of the moduli space)
is annihilated by a differential operator, coming from the restriction
property. Since $\SLE_{8/3}$ is revertible, so is the
conditional version; this gives a commutation condition in a
multiply connected domain. For $\SLE_2$, we know that the restriction construction is the scaling
limit of a revertible discrete model, hence the restriction weight
$h_1(x)=\alpha_2 B_x-\lambda_2 \Gamma_x$ should be
a solution of this functional equation. 
We now give a direct derivation of this fact.

\begin{Lem}
\begin{enumerate}
\item If there exists a non-vanishing function $\psi$ such that ${\mc M}_1\psi={\mc
M}_2\psi=0$, then $h_1,h_2$ satisfy:
$$(\check{\mc L}_1+2V'_x(y))(\alpha V'_y(x)+h_2)=(\check{\mc
L}_2+2V'_y(x))(\alpha V'_x(y)+h_1)$$
\item The following identity holds:
$$(\ell_x+2V'_x(y))\Gamma_y-\frac{V'''_x(y)}6
=(\ell_y+2V'_y(x))\Gamma_x-\frac{V'''_y(x)}6.
$$
\item The condition (i) is satisfied if $h_1(x,\dots)=h(x,\dots)$,
$h_2(y,\dots)=h(y,\dots)$, and
$$
h(x,\dots)=\alpha_\kappa B_x-\lambda_\kappa\Gamma_x+\ell_x f$$
where $f$ is a function on standard configurations that does not depend on $x,y$.
\end{enumerate}
\end{Lem}

\begin{proof}
(i)
If $\psi$ is such that ${\mc M}_1\psi={\mc
M}_2\psi=0$, then also ${\mc M}_0\psi=0$ where:
\begin{eqnarray*}
{\mc M}_0&=&[{\mc M}_1,{\mc M}_2]-2V'_y(x){\mc M_1}+2V'_x(y){\mc M_2}\\
&=&\left(\check{\mc L}_1 V_y(x)-\check{\mc
L}_2A_x+2V'_x(y)V_y(x)-2A_xV'_y(x)+\kappa\alpha V_y''(x)\right)\partial_x\\
&&-\left(\check{\mc L}_2 V_x(y)-\check{\mc
L}_1A_y+2V'_y(x)V_x(y)-2A_yV'_x(y)+\kappa\alpha
V_x''(y)\right)\partial_y\\
&&+\sum\left(\check{\mc L}_1 V_y(z_i)-\check{\mc L}_2 V_x(z_i)+2V'_x(y)V_y(z_i)-2V'_y(x)V_x(z_i)\right)\partial_{z_i}\\
&&+\left((\check{\mc L}_1+2V'_x(y))(\alpha V'_y(x)+h_2)-(\check{\mc
L}_2+2V'_y(x))(\alpha V'_x(y)+h_1)\right)
\end{eqnarray*}
The first-order terms vanish, from (\ref{cmu3}), (\ref{cmu4}).

(ii)
We want to interpret the left-hand side and the right-hand side of the
equation in terms of Brownian measures. Let us start with the
left-hand side. In the standard domain $D$, we grow a vertical slit
$\gamma$ at $x$. For small $t$, $H_t=D\setminus \gamma_{[0,t]}$ is
conformally equivalent to a standard domain $D_t$; the conformal
equivalence $f_t$ can be expanded in $t$ as:
$$f_t(z)= z+tV_x(z)+o(t)$$
where $V_x$ is the Schwarz kernel with pole at $x$ in the domain $D$. 
Hence $f'_t(y)^2=1+2tV_x'(y)+o(t)$ and $Sf_t(y)=tV'''_x(y)+o(t)$.
Let us consider
$\Gamma(H_t,y)$. From the covariance property of $\Gamma$, we get:
\begin{align*}
\Gamma(H_t,y)-\Gamma(H_0,y)&=f'_t(y)^2\Gamma(D_t,f_t(y))-\frac {Sf_t(y)}6-\Gamma(H_0,y)\\
&=t\left((V_x(y)\partial_y+\sum_i V_x(z_i)\partial_{z_i}+2V'_x(y))\Gamma_y-\frac{V'''_x(y)}6\right)+o(t)
\end{align*}
Now $\Gamma(H_t,y)-\Gamma(H_0,y)$ is the measure for Brownian bubbles
rooted at $y$ of bubbles that intersect $\gamma_{[0,t]}$ as well as
one of the holes. Consider $\tilde\gamma$ a vertical slit growing at
$y$, and the loop soup measure:
\begin{align*}
\mu^{loop}_\H\left(\{\delta:\delta\cap\gamma\neq\varnothing,\delta\cap\tilde\gamma\neq\varnothing,\delta\subsetneq
H_0\}\right)&=\mu^{bub}_x\left(\{\delta:\delta\cap\tilde\gamma\neq\varnothing,\delta\subsetneq
H_0\}\right)(t+o(t))\\
&=\mu^{bub}_y\left(\{\delta:\delta\cap\gamma\neq\varnothing,\delta\subsetneq
H_0\}\right)(t+o(t))
\end{align*}
as follows from Proposition 11 in \cite{LW}. Here we are considering
loops that intersect the two small slits $\gamma,\tilde\gamma$. We can
root such loops either near $x$ or $y$, which gives us the right-hand
side and the left hand-side of the claimed identity. 

(iii)
As we did earlier, we can expand the condition (i) at $y=x$, writing
$h_1=\alpha B_x-\lambda_\kappa \Gamma_x+h(x,\dots)$, $h_2=\alpha
B_y-\lambda_\kappa \Gamma_y +h(y,\dots)$ for some unknown function $h$.The
functional equation for $h$ is the linear equation :
$$(\ell_x+2V'_x(y))h(y)=
(\ell_y+2V'_y(x))h(x)$$
(from (\ref{cmu5}) and (ii)). Note that this equation does not involve
$\kappa$. Also, using the identity \ref{cmu3}, one can write:
$$(\ell_x+2V'_x(y))\ell_y+(\dots)\partial_x=
(\ell_y+2V'_y(x))\ell_x+(\dots)\partial_y.$$
So if $f$ does not depend on $x,y$, $h(x)=\ell_x f$ gives a solution
of the functional equation.
\end{proof}

Let us consider again the
equation ``without second member'':
$$(\ell_x+2V'_x(y))h(y)=(\ell_y+2V'_y(x))h(x)$$
appearing in (iii).
From (\ref{cmu3}),
we see that if $z$ is a marked point, then $h(x)=V_x(z)$  is a solution
for any marked point $z$. Also, differentiating the identity (\ref{cmu3}) w.r.t. $z$,
two terms $V'_x(z)V'_y(z)$ cancel out, and we get:
$$(\ell_x+2V'_x(y))V'_y(z)=(\ell_y+2V'_y(x))V'_x(z)$$
So $h(x)=V'_x(z)$ is a solution. Similarly,
$h(x)=(V_x(z_2)-V_x(z_1))/(z_2-z_1)$ is also a solution. This
corresponds to rational solutions in the simply connected case. It is
not so clear what would be $f$ in these cases.

As in Section 4, we can interpret a solution $h$ as the tangent map of
a cocycle. More precisely, if $h$ is as above, then one can define a
function $C$ on hulls and residual configurations (configurations
determined by the domain and the marked points $z_1,\dots$) such that:
\begin{enumerate}
\item For all hulls $A,B$,
 $C(A.B,z)=C(B,\phi_A(z))C(A,z)$.
\item If $A$ is a hull of half-plane capacity $\eps$ located at $x$, then
$C(A,z)=1+2\eps h(x,z)+o(\eps)$. 
\end{enumerate}
Some of the arguments we used in the simply connected case can be
replicated here; though a general classification of such cocycles
appears to be more difficult.

\vspace{1cm}

\noindent {\bf Acknowledgments.} I wish to thank Greg Lawler and
Wendelin Werner for stimulating and fruitful conversations.

\bibliographystyle{abbrv}
\bibliography{biblio}

\begin{thebibliography}{10}

\bibitem{AizBur}
M.~Aizenman and A.~Burchard.
\newblock H\"older regularity and dimension bounds for random curves.
\newblock {\em Duke Math. J.}, 99(3):419--453, 1999.

\bibitem{BauFr}
R.~O. Bauer and R.~M. Friedrich.
\newblock {On Chordal and Bilateral SLE in multiply connected domains}.

\bibitem{B1}
V.~Beffara.
\newblock { The dimension of SLE curves}.
\newblock {\em preprint, arXiv:math.PR/0211322}, 2002.

\bibitem{CamNew}
F.~Camia and C.~M. Newman.
\newblock {The Full Scaling Limit of Two-Dimensional Critical Percolation}.

\bibitem{Ca3c}
J.~Cardy.
\newblock Corrigendum: ``{S}tochastic {L}oewner evolution and {D}yson's
  circular ensembles'' [{J}. {P}hys. {A} {\bf 36} (2003), no. 24,
  {L}379--{L}386 ].
\newblock {\em J. Phys. A}, 36(49):12343, 2003.

\bibitem{Ca3}
J.~Cardy.
\newblock Stochastic {L}oewner evolution and {D}yson's circular ensembles.
\newblock {\em J. Phys. A}, 36(24):L379--L386, 2003.

\bibitem{D4}
J.~Dub{\'e}dat.
\newblock $\sle(\kappa,\rho)$ martingales and duality.
\newblock {\em Ann. Probab., to appear}, 2003.

\bibitem{D3}
J.~Dub{\'e}dat.
\newblock Critical percolation in annuli and {${\rm SLE}\sb 6$}.
\newblock {\em Comm. Math. Phys.}, 245(3):627--637, 2004.

\bibitem{D5}
J.~Dub{\'e}dat.
\newblock Excursion decompositions for $\sle$ and {W}atts' crossing formula.
\newblock {\em preprint, arXiv:math.PR/0405074}, 2004.

\bibitem{D7}
J.~Dub{\'e}dat.
\newblock Euler integrals for commuting {SLE}s.
\newblock {\em preprint, arXiv:math.PR/0507276}, 2005.

\bibitem{Fomin}
S.~Fomin.
\newblock Loop-erased walks and total positivity.
\newblock {\em Trans. Amer. Math. Soc.}, 353(9):3563--3583 (electronic), 2001.

\bibitem{KozL}
M.~J. Kozdron and G.~F. Lawler.
\newblock {Estimates of random walk exit probabilities and application to
  loop-erased random walk}.

\bibitem{LSW3}
G.~Lawler, O.~Schramm, and W.~Werner.
\newblock Conformal restriction: the chordal case.
\newblock {\em J. Amer. Math. Soc.}, 16(4):917--955 (electronic), 2003.

\bibitem{Law2}
G.~F. Lawler.
\newblock {The Laplacian-$b$ random walk and the Schramm-Loewner evolution }.

\bibitem{Law}
G.~F. Lawler.
\newblock {\em Conformally invariant processes in the plane}, volume 114 of
  {\em Mathematical Surveys and Monographs}.
\newblock American Mathematical Society, Providence, RI, 2005.

\bibitem{LSW2}
G.~F. Lawler, O.~Schramm, and W.~Werner.
\newblock Conformal invariance of planar loop-erased random walks and uniform
  spanning trees.
\newblock {\em Ann. Probab.}, 32(1B):939--995, 2004.

\bibitem{LW}
G.~F. Lawler and W.~Werner.
\newblock The {B}rownian loop soup.
\newblock {\em Probab. Theory Related Fields}, 128(4):565--588, 2004.

\bibitem{RS01}
S.~Rohde and O.~Schramm.
\newblock Basic properties of {SLE}.
\newblock {\em Ann. of Math. (2)}, 161(2):883--924, 2005.

\bibitem{S0}
O.~Schramm.
\newblock Scaling limits of loop-erased random walks and uniform spanning
  trees.
\newblock {\em Israel J. Math.}, 118:221--288, 2000.

\bibitem{S1}
O.~Schramm.
\newblock A percolation formula.
\newblock {\em Electron. Comm. Probab.}, 6:115--120 (electronic), 2001.

\bibitem{Sm1}
S.~Smirnov.
\newblock {Critical percolation in the plane. I. Conformal Invariance and
  Cardy's formula II. Continuum scaling limit}.
\newblock {\em in preparation}, 2001.

\bibitem{W3}
W.~Werner.
\newblock Conformal restriction and related questions.
\newblock {\em Lecture notes, ICMS Edinburgh, July 2003}, 2003.

\bibitem{W2}
W.~Werner.
\newblock Girsanov's transformation for {${\rm SLE}(\kappa,\rho)$} processes,
  intersection exponents and hiding exponents.
\newblock {\em Ann. Fac. Sci. Toulouse Math. (6)}, 13(1):121--147, 2004.

\bibitem{W1}
W.~Werner.
\newblock Random planar curves and {S}chramm-{L}oewner evolutions.
\newblock In {\em Lectures on probability theory and statistics}, volume 1840
  of {\em Lecture Notes in Math.}, pages 107--195. Springer, Berlin, 2004.

\bibitem{Wi}
D.~B. Wilson.
\newblock Generating random spanning trees more quickly than the cover time.
\newblock In {\em Proceedings of the Twenty-eighth Annual ACM Symposium on the
  Theory of Computing (Philadelphia, PA, 1996)}, pages 296--303, New York,
  1996. ACM.

\bibitem{DZ}
D.~Zhan.
\newblock {Stochastic Loewner Evolution in doubly connected domains}.
\newblock {\em Probab. Theory Related Fields}, 2003.

\bibitem{DZT}
D.~Zhan.
\newblock {Random Loewner chains in Riemann surfaces}, 2004.

\end{thebibliography}

-----------------------

Courant Institute

251 Mercer St., New York NY 10012

dubedat@cims.nyu.edu

\end{document}